%% file: aipsamp.tex
\documentclass{article}
\usepackage[a4paper, total={5.75in, 9in}]{geometry}
\usepackage{amsmath, amssymb}
\usepackage{graphicx}% Include figure files
\usepackage{dcolumn}% Align table columns on decimal point
\usepackage{bm} % bold math
\usepackage[utf8]{inputenc}
\usepackage[T1]{fontenc}
\usepackage{mathptmx}
\usepackage{mathrsfs}
\usepackage{etoolbox}
\usepackage{xcolor}

\usepackage{algorithm}
\usepackage{algorithmic}
\usepackage{colortbl}
\definecolor{bluePoli}{cmyk}{0.4,0.1,0,0.4}

\newcommand{\pspace}{\Theta}  % spazio dei parametri
\newcommand{\nmu}{p} % # parametri
\newcommand{\mub}{\boldsymbol{\mu}} % vettore dei parametri
\newcommand{\encoder}{\mathcal{E}}
\newcommand{\decoder}{\mathcal{D}}
\newcommand{\processor}{\mathcal{P}}
\newcommand{\vb}{\boldsymbol{v}}
\newcommand{\eb}{\boldsymbol{e}}
\newcommand{\ub}{\boldsymbol{u}}
\newcommand{\xib}{\mathbf{\xi}}
\newcommand{\x}{\mathbf{x}}
\newcommand{\Ee}{\encoder_{e}}
\newcommand{\Ev}{\encoder_{v}}
\newcommand{\ntrain}{N_{\text{train}}}
\newcommand{\nfom}{N_{\mub}^{h}}
\newcommand{\mesh}{\mathcal{M}_{\mub}^{h}}
\newcommand{\meshi}{\mathcal{M}_{\mub_{i}}^{h}}
\newcommand{\nicola}[1]{#1}

\makeatletter
\def\@email#1#2{%
 \endgroup
 \patchcmd{\titleblock@produce}
  {\frontmatter@RRAPformat}
  {\frontmatter@RRAPformat{\produce@RRAP{*#1\href{mailto:#2}{#2}}}\frontmatter@RRAPformat}
  {}{}
}%
\makeatother

\begin{document}

%\preprint{AIP/123-QED}

\title{Deep Learning-based surrogate models for\\parametrized PDEs: 
%including geometrical features 
handling geometric variability\\through graph neural networks}

\author{Nicola Rares Franco$^{1}$, Stefania Fresca$^{1}$, Filippo Tombari$^{1}$, Andrea Manzoni$^{1}$}



\date{$^{1}$ MOX, Department of Mathematics, Politecnico di Milano, Milan, Italy}

\maketitle

\begin{abstract}
Mesh-based simulations play a key role when modeling complex physical systems that, in
many disciplines across science and engineering, require the solution of parametrized time-dependent nonlinear partial differential equations (PDEs). In this context, full order models (FOMs), such as those relying on the finite element method, can reach high levels of accuracy, however often yielding intensive simulations to run. For this reason, surrogate models are developed to replace computationally expensive solvers with more efficient ones, which can strike favorable trade-offs between accuracy and efficiency. 
This work explores the potential usage of graph neural networks (GNNs) for the simulation of time-dependent PDEs in the presence of 
geometrical variability. In particular, we propose a systematic strategy to build surrogate models based on a data-driven time-stepping scheme where a GNN architecture is used to efficiently evolve the system. 
With respect to the majority of surrogate models, the proposed approach stands out for its ability of tackling problems with parameter dependent spatial domains, while simultaneously generalizing to different geometries and mesh resolutions.
We assess the effectiveness of the proposed approach through a series of numerical experiments, involving both two- and three-dimensional problems, showing that GNNs can provide a valid alternative to traditional surrogate models in terms of computational efficiency and generalization to new scenarios. We also assess, from a numerical standpoint, the importance of using GNNs, rather than classical dense deep neural networks, for the proposed framework.
\end{abstract}

\;\\\noindent\textbf{Short summary}\;\;\;\;\;
Geometric variability is a major obstacle in surrogate modeling, as classical approaches, such as the reduced basis method, can account for such degree of complexity only under severe simplifications, in an intrusive way, featuring remarkable computational costs. In this paper, we propose the use of graph neural networks to efficiently evolve dynamical systems defined on different domains and geometries. The networks are trained on a collection of trusted samples, obtained through accurate numerical simulations, and are shown to be capable of generalizing to unseen geometries without loss of accuracy. Despite assessed on a series of simplified test cases, numerical results suggest that the proposed approach can pave a new way for handling geometric variability in surrogate modeling, potentially leading to novel methodologies capable of combining GNNs and classical techniques.

\section{\label{sec:intro}Introduction} 

Thanks to accurate and reliable numerical simulations, we are now able to simulate, monitor and forecast very complex physical phenomena such as those arising in computational physics, biology and engineering.
However, when it comes to many-query applications, such as, e.g., optimal control and uncertainty quantification tasks, the elevated computational cost constitutes a major limitation that hinders the effective potential of numerical simulations.

As already explored by several researchers, one way to overcome this complexity is to rely on \textit{surrogate models}: suitable emulators that are capable of replicating the outputs of classical PDE solvers - thereby referred to as Full Order Models (FOM) - at a reduced computational cost. This practice is also known as Reduced Order Modeling (ROM). As of today, domain practitioners can count on a very large number of ROM techniques, each with its own advantages and limitations. Just to mention some of them, these include: intrusive and non-intrusive projection-based ROMs \cite{negri2013reduced, quarteroni2015reduced, hesthaven2016certified, hesthaven2018non, guo2019data, guo2022bayesian}, which can effectively tackle diffusive problems, especially in the case of affinely parametrized operators; adaptive methods based on, e.g., ROM augmentation \cite{carlberg2015adaptive, hesthaven2022rank}, clustering \cite{hess2019localized}, interpolation \cite{amsallem2011online} or space-time splittings \cite{pagliantini2021dynamical}, which are particularly suited for modeling shock waves, Hamiltonian systems, etc.; nonlinear reduction techniques based on, e.g., spectral submanifolds \cite{buza2021using, li2023model} and library representations \cite{bonito2021nonlinear}, which can provide users with solid theoretical guarantees; Deep-Learning based ROMs (DL-ROMs) relying on deep autoencoders, which, if provided with enough data, can address both stationary and time-dependent problems, even in the presence of severe nonlinearities and singular behaviors \cite{fresca2021comprehensive, franco2023deep, fresca2022pod, brivio2023error, cicci2022deep, fatone2022long, romor2023non}.
\\\\
All these approaches are grounded on a common assumption, that is: the underlying FOM must be identified once and for all, with a fixed spatial discretization and a precise number of degrees of freedom (dofs) $N_{h}$.
This fact, however, poses a major limitation when having to deal with PDEs defined over parametrized domains. 

Assume for instance that the governing equations depend on a vector of geometrical parameters $\mub$, which can affect the shape and the configuration of the underlying spatial domain $\Omega=\Omega_{\mub}.$ Then, at the discrete level, each $\mub$ instance will correspond to a suitable high-fidelity mesh $\mesh$ entailing $\nfom$ dofs. The issue, here, is that as soon as we change the value of the geometrical parameters, say from $\mub$ to $\mub'\neq\mub$, the total number of dofs might change, $N_{\mub'}^{h}\neq\nfom$; furthermore, in general, even if $N_{\mub'}^{h}=\nfom$, we will not be able to match the dofs in the two meshes. For projection-based ROMs, this makes the construction of a unique projection matrix $\mathbf{V}\in\mathbb{R}^{N_{h}\times n}$ impossible; similarly, we cannot rely on naive DL-ROMs as these would require the construction of an autoencoder network $\Psi\circ\Psi'$ with $\Psi':\mathbb{R}^{N_{h}}\to\mathbb{R}^{n}$ and $\Psi':\mathbb{R}^{n}\to\mathbb{R}^{N_{h}}$. Local techniques based on clustering algorithms can partially resolve this issue by providing a different projector $\mathbf{V}_{i}$ for each parametric instance $\mub_{i}$ observed in the so-called \textit{offline stage}, with $i=1,\dots,q$. However, this approach would incur in severe limitations during its \textit{online} usage, as the resulting ROM would not be applicable whenever a new parametric instance $\mub\notin\{\mub_{i}\}_{i=1}^{q}$ is given. Similar issues are encountered when dealing with FOMs that use mesh-adaptive strategies, even if the geometry is kept fixed.

The purpose of this work is to overcome these limitations by relying on Graph Neural Networks (GNNs), thus providing a flexible approach to surrogate modeling that is capable of handling geometric variability and generalizing to unseen geometries. The idea is inspired by the recent successes of GNNs in scientific applications \cite{shlomi2020graph, shukla2022scalable, horie2022physics} and shares some similarities, that we discuss below, with other recent works.

\subsection{GNNs in surrogate and reduced order modeling}

GNNs are a particular class of neural network architectures that were originally proposed as a way to handle statistical data defined over graphs \cite{scarselli2008graph, battaglia2018relational}. To simplify, given a (directed) graph $G=(V,E)$ with vertices $V$ and edges $E\subseteq V\times V$, a GNN is a computational unit that can receive a set of node features at input, $\vb:V\to\mathbb{R}^{l}$, and return a corresponding set of node features at output $\vb':V\to\mathbb{R}^{l'}.$ The same GNN unit can process data coming from different graphs. The only restrictions are: i) the size of the input and output features, $l$ and $l'$, respectively; ii) the fact that each input-output pair must be defined over the same graph.

This is possible because, differently from other architectures such as dense deep feed forward networks (DNNs), GNNs adopt a local perspective: information is processed at the nodal level through a combination of message-passing steps (communication of nearby nodes) and aggregation routines. This added flexibility makes GNN capable of handling data defined over different graphs and, eventually, provides them with the ability to generalize over unseen geometries. In the Deep-Learning literature, this fact is known as \textit{relational inductive bias}. In general, the term \textit{inductive bias} refers to the ability of a learning algorithm to prioritize one solution (or interpretation) over another, independently of the observed training data, and it can express (explicitly or implicitly) assumptions about either the data-generating process or the space
of solutions\cite{battaglia2018relational}. In the case of GNNs, the implicit assumption is that the output of a given neuron is primarily affected by its neighbouring neurons (thus the term \textit{relational}), so that local effects are stronger than global ones.
\\\\
Our idea for the present work is to exploit the capabilities of GNNs in order to learn a nonintrusive data-driven time-stepping scheme for evolving high dimensional parameter dependent dynamical systems. To this end, we interpret discrete FOM solutions $$\mathbf{u}_{\mub}=\left[\mathrm{u}_{\mub,1},\dots,\mathrm{u}_{\mub,\nfom}\right]^{T}\in\mathbb{R}^{\nfom}$$ as collections of nodal features $\ub_{\mub}:V_{\mub}\to\mathbb{R}$, were $$V_{\mub}=\left\{\x_{\mub,i}\right\}_{i=1}^{\nfom}$$ are the vertices of the underlying mesh (sorted coherently with the FOM dofs), so that
$$\ub_{\mub}\left(\x_{\mub,i}\right):=\mathrm{u}_{\mub}^{(i)}.$$
Then, this graph-mesh equivalence allows us to construct a GNN module that can evolve discrete solutions defined over different meshes (and different domains).
\\\\
\nicola{This work} %Our approach 
finds its main inspiration in a recent contribution by Pfaff et al.\cite{pfaff2020learning}, where the authors propose a GNN architecture for learning mesh-based simulations in a time-dependent framework. 
\nicola{Our purpose is to transpose their ideas to the realm of ROM for parametrized PDEs, and to propose a systematic approach for handling geometric variability. To this end, we shall adopt a purely mathematical perspective, as to convey the overall idea in the language that ROM practitioners are mostly familiar with.} 

\nicola{Nonetheless, aside from the surrounding framework and the mathematical formalism, our proposal also features a few practical differences with respect to the work by Pffaf et al., namely: i) the introduction of \textit{global features}, which we use to extend the overall approach  to nonautonomous systems and, possibly, to PDEs that depend both on physical and geometrical parameters; ii) the definition of the loss function,  which we complement with an additional term concerning the approximation of the time-derivative; iii) an explicit superimposition of a Runge-Kutta-like time-stepping scheme.}

%The main difference with respect to our proposal lies in the surrounding framework, as our work systematically addresses the case of parametrized systems. Indeed, despite this aspect is not explicitly addressed in the present work and is left to forthcoming publications, the proposed strategy is designed to handle nonlinear PDEs depending on both physical and geometrical parameters. Furthermore, the two approaches also have a few practical differences, namely:  and the explicit superimposition of a Runge-Kutta-like time-stepping scheme.

In this sense, our work is much closer to the one by Pegolotti et al.\cite{pegolotti2023learning}, where the authors explore the use of GNNs for reduced order modeling of cardiovascular systems. Still, their framework remains quite different from ours as they only consider a fixed number of possible geometries, thus not allowing for a continuous parametrization, and they focus on a specific physical system.
A more flexible use of GNNs is found in the recent contribution by Gladstone et al.\cite{gladstone2023gnn}, where a similar paradigms is exploited to surrogate classical PDE solvers. Their analysis, however, is limited to time-independent PDEs and does not transfer to dynamical systems.

Finally, for what concerns surrogate and reduced order modeling, we mention that some authors are also exploring the integration of GNNs together with ROM techniques: see, e.g., the GCA-ROM, a GNN-variation of the DL-ROM approach recently proposed by Pichi et al.\cite{pichi2023graph} Nonetheless, these techniques are extremely different with respect to our proposal, as, in order to tackle both stationary and time-dependent PDEs, they neglect the dynamical nature of the system, that is: they treat time as an additional parameter, thus ignoring the Markovian structure that characterizes the majority of evolution equations.

\subsection{Outline of the paper}
The paper is organized as follows. First, in Section \ref{sec:setup}, we formally introduce the problem of surrogate modeling for parametrized dynamical systems. Then, in Section \ref{sec:gnns}, we provide the reader with the fundamental building blocks required for our construction and present the corresponding GNN architectures. We then put things into action in Section \ref{sec:surrogate}, where we dive into the details of the proposed approach. Finally, we devote Section \ref{sec:exp} to the numerical experiments.

\section{\label{sec:setup}Modeling time-dependent PDEs}

%For the sake of generality, despite the numerical experiments will refer to linear time-dependent parametrized PDEs, we here consider the nonlinear case, since the approach presented in the following can be easily extended to more general problems. More precisely, 
We consider a PDE system depending on a set of
input parameters $ \mub \in \pspace $, where the parameter space $ \pspace \subset \mathbb{R}^{\nmu}$ is a bounded and closed set; in our analysis, input parameters may represent both physical and geometrical properties of the system, like, e.g., material properties, boundary conditions, or the shape of the domain itself. For the time being, however, we focus on the treatment of geometrical parameters, since the extension to the case where both physical and geometrical parameters is straightforward.
Throughout the paper, we adopt a fully algebraic perspective and assume that the governing equations have already been discretized in space by means of a suitable high-fidelity approximation -- which, here, is allowed to depend on $\mub$ -- such as, e.g., the finite element method. 
%to start fromthe high-fidelity (spatial) approximation of the PDE system considered.
Regardless of the spatial discretization adopted, the FOM can be expressed as a nonlinear, high-dimensional 
parametrized dynamical system. Hence, given $\mub \in \pspace$, we aim at solving the initial value problem:
\begin{equation}
    \label{eq:pde}
    \left\{
    \begin{aligned}
    &\mathbf{\dot u}_{\mub}(t)  = \mathbf{f}\left(t, \mathbf{u}_{\mub}(t), \mub\right), \qquad t \in (0,T) , \\
    &\mathbf{u}_{\mub}(0)  = \mathbf{g}_{\mub},
    \end{aligned}
    \right.
\end{equation}
where $\mathbf{u}_{\mub}: [0,T) \rightarrow \mathbb{R}^{\nfom}$ is the parametric solution to (\ref{eq:pde}), while $$\mathbf{g}_{\mub} \in \mathbb{R}^{\nfom}\quad\text{and}\quad\mathbf{f}(\cdot,\cdot,\mub): (0,T) \times \mathbb{R}^{\nfom}  \rightarrow \mathbb{R}^{\nfom}$$ are the initial condition and a - possibly nonlinear - function encoding the dynamics of the system, respectively. The FOM dimension, $\nfom$, is related to the finite dimensional subspaces introduced for the sake of space discretization – here $h > 0$ denotes a discretization parameter, such as the maximum diameter of the elements
in the computational mesh $\mesh$% related to a triangulation $\mathcal{T}_h(\mub)$ such that $\Omega_h(\mub) = \textnormal{int}(\cup_{K \in \mathcal{T}_h(\mub )}K)$ù
; consequently, $\nfom$ can be extremely large if the PDE problem describes complex physical behaviors and/or high degrees of accuracy are required for its solution. Furthermore, the number of degrees of freedom (dofs) of the problem may depend on the geometrical parameters contained in $\mub$ since, by modifying their values, the number of vertices in the computational mesh $\mesh$ can vary. %Consequently, we strive to discover solutions that can generalize across meshes with varying node counts.
We thus aim at approximating the set
\begin{equation}
\label{eq:manifold}
    \mathcal{S} = \{ \mathbf{u}_{\mub}(t) | \ t \in [0,T), \ \mub \in \pspace \subset \mathbb{R}^{p} \} \subset \bigcup_{\mub \in \pspace}\mathbb{R}^{\nfom}
\end{equation}
of the solutions to (\ref{eq:pde}) when $(t; \mub)$ varies in $[0, T) \times \pspace$, also referred to as {solution manifold}.

In order to numerically approximate problem (\ref{eq:pde}), even at the FOM level, one must rely on suitable time-integration schemes, such as the backward differentation formulae \cite{quarteroni1994numerical}. Thus, having fixed a uniform partition of $(0,T)$ in $N_{t}$ equally spaced subintervals, and by denoting with $\mathbf{u}^n$, the solution $\mathbf{u}$ at time $t^n = n\Delta t$, where $\Delta t:=T/N_{t}$, our ultimate aim is to solve:
\begin{equation}
\label{eq:fem}
 \left\{
    \begin{aligned}
    & \frac{\mathbf{u}^{n+1}_{\mub} - \mathbf{u}^{n}_{\mub}}{\Delta t}   = \mathbf{f}\left(t^{n+1}, \mathbf{ u}^{n+1}_{\mub}, \mub\right), \qquad n \geq 0, \\
    &\mathbf{u}^{0}_{\mub}  = \mathbf{g}_{\mub}.
    \end{aligned}
    \right.
\end{equation}
Equation (\ref{eq:fem}) requires the solution, at each time instance, of a  nonlinear system depending on the input parameter vector $\mub$ which may entail high computational times, especially when dealing with a multi-query or real-time context. To achieve computational efficiency, multi-query analysis and real-time problems must rely on suitable surrogate models which can be built according to different strategies. Motivated by this, our goal is the efficient approximation of the solution manifold in (\ref{eq:manifold}) by decreasing the complexity related to the solution of the FOM and preserving high level of accuracy.

In this work, we introduce a Deep Learning-based surrogate model that exploits graph neural networks (GNNs) \cite{hamilton2017inductive} to efficiently evolve the time-discrete dynamical system in \eqref{eq:fem}. Here, the use of GNNs is motivated by their unique ability of handling data defined on different graphs/meshes, which can result in extremely flexible models capable of generalizing to new, unseen, geometries and spatial resolutions (in the Deep Learning literature, this fact is usually referred to as relational inductive bias\cite{battaglia2018relational}).
%GNNs allow to find an approximate solution $\Tilde{\mathbf{u}}(t,\mub) \approx \mathbf{u}(t,\mub)$ which is not strictly dependent on the space discretization parameter $h$, or the triangulation $\mathcal{T}_h$, related to the computational mesh. This choice is driven by the  of graph-based DL algorithms making the model proposed capable to handle geometric variability in the solution of problem \ref{eq:fem} and thus generalize to new, unseen geometries. More precisely, we introduce 
In mathematical terms, we aim at constructing a GNN architecture $\Phi$ %$\bm{\Phi}: \bigcup_{\mub \in \mathcal{P}}\mathbb{R}^{\nfom} \times \mathcal{P} \rightarrow \bigcup_{\mub \in \mathcal{P}}\mathbb{R}^{\nfom} $ %\Tilde{\mathcal{S}}_h,$
such that:
$$\mathbf{u}^{n+1}_{\mub}  \approx \bm{\Phi}(\mathbf{u}^{n}_{\mub}, t^{n}, \mub).$$
\nicola{From an abstract point of view, the above can be seen as an extension of the MeshGraphNet model as originally proposed by Pfaff et al.\cite{pfaff2020learning}: here, in fact, the time variable is included explicitly, which allows us to address the more general case of nonautonomous systems.}

%\begin{equation}
% \left\{
%    \begin{aligned}
%    &\Tilde{\mathbf{u}}^{n+1}_{\mub}   = \bm{\Phi}(\Tilde{\mathbf{u}}^{n}_{\mub}; \mub) \quad n = 0,\ldots,N - 1, \\
%    &\Tilde{\mathbf{u}}^0(\mub)  = \mathbf{u}^0(\mub).
%    \end{aligned}
%    \right.
%\end{equation}

\section{\label{sec:gnns}Graph Neural Networks}
GNNs were initially conceived as an extension of convolutional neural networks (CNNs) to operate on graph-structured data and overcome their limitations in this domain. In graph theory, graphs are used to describe systems made of nodes and their connections (edges). An image can be regarded as a graph with regular and well-organized connections in the Euclidean space, where each pixel corresponds to a node in the graph. In this particular case, the aforementioned CNN architectures can exploit the
peculiar structure of the graph to extract meaningful spatial features. However, these models become inapplicable as soon as the structure of the underlying graph becomes slightly more sophisticated: in practice, this fact has remarkable consequences as in many real-time applications (e.g., traffic networks, socialnetworks) data are naturally defined over general, possibly non-Euclidean, graphs.

\begin{figure}
    \centering
  \includegraphics[width=0.9\textwidth]{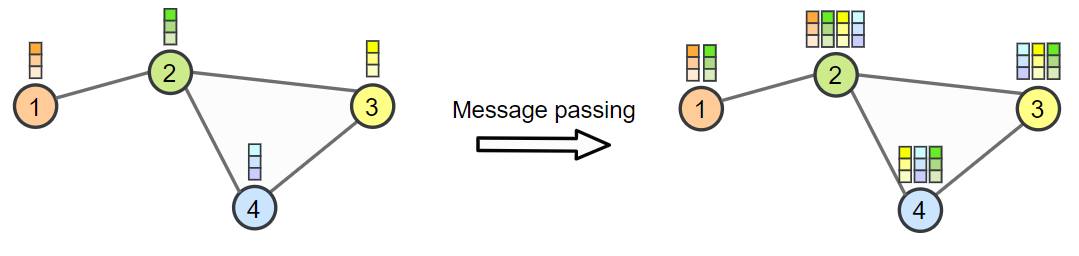}
  \caption{\small{Message propagation and aggregation. The information is broadcasted from different levels of depth of the graph. For each node, at each message passing step,the information is collected from the neighbors and aggregated. In this way, adding message-passing steps can be seen as connecting nodes which can also be far from each other. Figure courtesy of Phillip Lippe (University of Amsterdam, QUVA lab).} }
  \label{fig:mp}
\end{figure}

To promote the use of Deep Learning in those applications, GNNs were developed to extract spatial features over general graphs, by inspecting neighboring nodes, with arbitrary connection in a non-Euclidean space \cite{pmlr-v48-niepert16,zhang2019graph}. Thus, GNNs can be considered as a generalized version of CNNs over general graphs. Clearly, this generalization comes with some differences between the two architectures. For instance, while the outputs of CNNs are affected by the ordering of the pixels, the same is not true for GNNs, as the action of the latter is uniquely determined by the connectivity of the graph (in this concern, note, for instance, that the connectivity of an image remains the same even if we flip it either vertically or horizontally). Moreover, GNNs adopt a \textit{graph-in, graph-out} architecture meaning that these model types accept a graph as input, with information loaded into its nodes and edges, and progressively transform these embeddings, without changing the connectivity of the input graph: in contrast, CNN layers usually modify the resolution of the input upon their action.
\\\\
The fundamental ingredient of a basic GNN layer is the so-called \textit{message passing} operation, which enables the aggregation of node information while leveraging the depth of the graph. More precisely, a message passing step consists of two components:
\begin{itemize}
    \item message computation: each node creates a message to be sent to other nodes later;
    \item aggregation: each node aggregates the messages from the neighborhood. 
\end{itemize}
This message-passing propagation can be seen as an information retrieval task from different levels of depth of the graph. A simple visualization of the message propagation is shown in Figure \ref{fig:mp}. For each node, the information comes from the neighbors. In this way, adding message-passing steps can be seen as connecting nodes that can be also far from each other. 
Graph-based algorithms, such as, e.g, graph convolutional networks \cite{kipf2016semi}, GraphSage \cite{HamiltonYL17}, graph attention networks \cite{veli2017graph}, graph transformer operator \cite{ijcai2021p214} and interaction networks \cite{10.5555/3157382.3157601}, differ in the way the message is computed and the aggregation is performed. In particular, depending on the chosen framework, the message-passing step may also involve the edges of the graph, where a corresponding set of \textit{edge features} can be loaded: in the next few pages, we shall describe this situation in full mathematical detail, as it will be of key importance for our construction.

In order to perform a message-passing operation, GNNs leverage on suitable data structures for representing the topology and connectivity of the graph. %it is necessary to find a suitable representation of the graph structure, and, especially, of its connectivity. %Graphs have different types of information that shall potentially be exploited to make predictions, i.e, nodes, edges, and connectivity. The first two features are relatively straightforward: for example, by means of nodes, we can form a node feature matrix $\mathbf{N}$ by assigning each node an index $i$ and storing the feature for node $i$ in $\mathbf{N}$. While these matrices can show a quite different nature, they can be processed without any special technique.\\
%However, representing graph connectivity is more complicated. 
In this concern, a classical choice is to exploit the edge connectivity matrix, that is, a $n_{edges} \times 2$ matrix where each row $k$ contains the indices of the source and destination nodes of the $k$th edge; this allows GNNs to stash the overall topology of the graph with a memory complexity of $\mathcal{O}(n_{edges})$. Roughly speaking, this is equivalent to storing a sparse version of the adjacency matrix of the graph, which, in principle, consists of $\mathcal{O}(n_{nodes}^2)$ entries.
\\\\
Before coming to our own use of GNNs for surrogate modeling, within this Section we take the chance to present some of the fundamental ingredients required for our construction. In particular, we shall describe in mathematical terms the concept of message-passing, and we shall introduce a particular GNN architecture known as the Encoder-Processor-Decoder model.

\subsection{\label{subsec:mp} The message-passing block: formal definition}

Given $l\in\mathbb{N}$, a graph-forward-pass with $l$ hidden features is a computational unit $F=F(\vb,\eb,G)$ that takes as input
\begin{itemize}
    \item [i)] a directed graph structure, $G=(V,E)$;
    \item [ii)] a collection of vertex features, $\vb:V\to\mathbb{R}^{l}$;
    \item [iii)] a collection of edge features, $\eb:E\to\mathbb{R}^{l}$ 
\end{itemize}
and outputs a new collection of vertex features with $l$-features per node, namely
$$F(\vb,\eb,G):V\to\mathbb{R}^{l}.$$
We think of $F$ as an object that transforms the vertex features associated to the nodes in the graph.

In GNN architectures, a message-passing block is a particular type of graph-forward-pass routine that exploits the local structure of the input graph $G$, only allowing communication of nearby nodes. Specifically, a message-passing block $F$ is comprised of two Multi-Layer Perceptron\cite{murtagh1991multilayer} (MLP) units, 
$$\psi_{v}:\mathbb{R}^{2l}\to\mathbb{R}^{l}\quad\text{and}\quad\psi_{e}:\mathbb{R}^{3l}\to\mathbb{R}^{l},$$
that completely characterize the action of $F$. However, in order to properly explain how the forward pass is carried out, we first need to introduce some notation. Given a graph $G=(V,E)$ and a collection of vertex features $\vb:V\to\mathbb{R}^{l}$, we write $\vb_{\text{in}}$ and $\vb_{\text{out}}$ for the maps
$$\vb_{\text{in}}:E\to\mathbb{R}^{l}\quad\quad \vb_{\text{out}}:E\to\mathbb{R}^{l}$$
given by
$$\vb_{\text{in}}(v_{1},v_{2}):=\vb(v_{1}),\quad \vb_{\text{out}}(v_{1},v_{2}):=\vb(v_{2}),$$
respectively, where $(v_{1},v_{2})\in E$ represents an oriented edge going from $v_{1}$ to $v_{2}$. In other words, passing from $\vb$ to $\vb_{\text{in}}$ is equivalent to transferring the information from the nodes to the edges, with the convention that a given edge inherits the features from its own \textit{source node}. Similarly, going from $\vb$ to $\vb_{\text{out}}$ is a way for storing the information about the \textit{destination nodes}. 

In the same spirit, it is also useful to define the dual operation, which transfers information from the edges to the nodes. More precisely, given $\eb:E\to\mathbb{R}^{l}$, we shall write $\overline{\eb}$ for the map $\overline{\eb}:V\to\mathbb{R}^{l}$
defined as
$$\overline{\eb}(v):=\sum_{(v_{1},v)\in E}\eb(v_{1},v),$$
that is, to go from $\eb$ to $\overline{\eb}$, we collapse all the features corresponding to edges with the same destination node.
 
Lastly, we shall denote by $\oplus$ the concatenation operator. Specifically, given any two functions with a common domain, e.g., $\boldsymbol{f}:X\to\mathbb{R}^{a}$ and $\boldsymbol{g}:X\to\mathbb{R}^{b}$, we write $\boldsymbol{f}\oplus\boldsymbol{g}$ to intend the map from $\Omega$ to $\mathbb{R}^{a+b}$ given by
$$\boldsymbol{f}\oplus\boldsymbol{g}(x):=[f_{1}(x),\dots,f_{a}(x),g_{1}(x),\dots,g_{b}(x)],$$
where $x\in X$ is a generic input, while $f_{i}$ and $g_{j}$ are the $i$th and $j$th components at the output of $\boldsymbol{f}$ and $\boldsymbol{g}$, respectively.
\\\\
We now have all the ingredients to rigorously define the forward-pass of a message-passing block. The action of a message-passing block $F$ with $l$ hidden features and computational units $\psi_{v},\psi_{e}$, is defined as
\begin{equation}
\label{eq:message-passing}
F(\vb,\eb,G)=\psi_{v}\circ\left(\vb\oplus\overline{\psi_{e}\circ(\eb\oplus \vb_{\text{in}} \oplus\vb_{\text{out}})}\right),
\end{equation}
where, as usual, $\circ$ denotes functional composition.

In plain words, Eq. \eqref{eq:message-passing} states that the vertex features at output, $F(\vb,\eb,G)$, are obtained as follows: first, the information available in the graph vertices is transferred to the edges and concatenated with the existing features, $\eb\oplus \vb_{\text{in}} \oplus\vb_{\text{out}}$; then, an MLP, $\psi_{e}$, is applied to the extended features to extract meaningful information; the latter, is then transferred back to the node vertices, yielding $\overline{\psi_{e}\circ(\eb\oplus \vb_{\text{in}} \oplus\vb_{\text{out}})}$. These hidden features -- which now live of the graph vertices -- are then appended to the original ones and later fed to a terminal MLP block, here given by $\psi_{v}$.

In general, we remark that the action of a message-passing step $(\vb,\eb,G)\mapsto F(\vb,\eb,G)$ is nonlinear because of the two MLPs, $\psi_{v}$ and $\psi_{e}$, entering the pipeline.
\newtheorem{remark}{Remark}
\begin{remark}
    Note that the operations $\vb\mapsto \vb_{\text{in}}$, $\vb\mapsto \vb_{\text{out}}$ and $\eb\mapsto \overline{\eb}$, require an exact knowledge of the graph structure $G$. Here, this fact is left implicit to keep the notation lighter.
\end{remark}

\begin{remark}
    In the literature, GNNs have been defined in several ways. One major difference lies in that some authors only talk about «node features», without contemplating the existence of «edge features». Here, we are adopting one of the most recent formulation of GNNs, as proposed by Battaglia et al.\cite{battaglia2018relational} Nonetheless, we believe that finding connections between different definitions can enhance practical understanding. To this end, we mention that in the classical formulation by Scarselli et al. (that is, without edge features), a major role is played by the «aggregation step», in which information coming from neighbouring nodes is collapsed onto a single value, e.g. via summation (see Equation 3 in the original work by Scarselli et al. \cite{scarselli2008graph}). Here, the same effect can be obtained via $\vb\mapsto\overline{\vb_{\text{in}}\oplus\vb_{\text{out}}}-(\boldsymbol{n}-1)\cdot\vb$, where $\boldsymbol{n}:V\to\mathbb{N}$ is a feature map that returns the connectivity of each node, whereas $\cdot$ stands for pairwise multiplication.
\end{remark}

\subsection{\label{subsec:model}The Encoder-Processor-Decoder model}
\label{sec:epd}
The Encoder-Processor-Decoder model is a powerful GNN-based architecture that can process mesh-based data \cite{battaglia2018relational,pfaff2020learning,sanchez2020learning}. More precisely, the latter accepts as input:
\begin{itemize}
    \item[i)] a directed graph $G=(V,E)$ associated to some mesh $\mathscr{M}$ embedded in a suitable ambient space $\mathbb{R}^{d}$, so that $V\subset\mathbb{R}^{d}$;
    \item[ii)] an input signal defined over the mesh vertices, namely $\boldsymbol{u}:V\to\mathbb{R}^{q}$;
    \item[iii)] a global feature vector, $\xib\in\mathbb{R}^{s}$, describing a given nonspatial property of the system (e.g., time).
\end{itemize}
Then, the output of such a model is a new signal $\boldsymbol{u}':V\to\mathbb{R}^{q'}$ defined over the given mesh.

As the name suggests, the Encoder-Processor-Decoder model is comprised of three modules, which we explain in detail below. These are all characterized by a common hidden-dimension, $l\in\mathbb{N}$, which we assume to be fixed hereon.

\subsubsection{Encoder module}
The encoder module is used to preprocess the input data and return a collection of hidden features defined, respectively, over the graph vertices $\Ev=\Ev(\ub,\xib,G)$ and the edge vertices $\Ee=\Ee(G)$. The two are obtained as follows.
\\\\
The node features $\Ev(\ub,\xib,G)$, are computed by combining a fixed nonlearnable transformation together with an MLP unit $\Psi_{\encoder}^{v}:\mathbb{R}^{q+s+1}\to\mathbb{R}^{l}$ that maps onto the hidden-state space. The former  has the purpose of expanding the node features with information coming from the global variables, $\xib$, and the graph $G$. More precisely, let $\boldsymbol{b}_{G}:V\to \{0,1\}$ be a flag for those nodes that lie on the boundary of the mesh, i.e.,  $\boldsymbol{b}_{G}(v)=1$ if and only if $v$ is a boundary vertex. Then, the action of $\Ev$ reads
\begin{equation}
\label{eq:encoder1}\Ev(\ub,\xib,G):=\Psi_{\encoder}^{v}\circ(\ub\oplus\xib\oplus\boldsymbol{b}_{G}),
\end{equation}
so that $\Ev(\ub,\xib,G):V\to\mathbb{R}^{l}.$ 
Here, with little abuse of notation, we have identified the vector $\xib$ with a constant map defined over $V$. As we mentioned, the preliminary transformation $\ub\mapsto \ub\oplus\xib\oplus\boldsymbol{b}_{G}$ is nonlearnable and has the sole purpose of augmenting the nodal features; conversely, the MLP unit introduces a learnable block that is optimized during training.

The edge features $\Ee(G)$ are computed following similar ideas. First, a set of nonlearnable features $e_{G}:E\to\mathbb{R}^{d+1}$ is extracted starting from the mesh coordinates. This is achieved by letting
$$e_{G}(\x_{1},\x_{2}):=\left[\frac{x_{1}^{(1)}+x_{2}^{1}}{2},\dots,\frac{x_{1}^{(d)}+x_{2}^{d}}{2},|\x_{1}-\x_{2}|\right]$$
where $(\x_{1},\x_{2})\in E.$ In other words, $e_{G}$ maps each edge to a vector containing the coordinates of its midpoint together with the edge length. These preliminary features are then fed to an MLP $\Psi_{\encoder}^{e}:\mathbb{R}^{d+1}\to\mathbb{R}^{l}$, i.e.
\begin{equation}
    \label{eq:encoder2}
    \Ee(\ub,\xib,G):=\Psi_{\encoder}^{e}\circ e_{G},
\end{equation}
which returns the encoded edge features.

\subsubsection{Processor module}
\label{sec:processor}
The encoded features, $\vb:=\Ev(\ub,\xib,G)$ and $\eb:=\Ee(G)$, are then elaborated by a GNN-based unit, called the processor $\processor$. The latter consists of $m$ message-passing-blocks, $F_{1},\dots,F_{m}$, each one acting as in \eqref{eq:message-passing}.  More precisely, the output of the processor module is given by
\begin{equation}
    \label{eq:processor}
    \processor(\vb,\eb,G):= F_{m}(\boldsymbol{h}_{m},\eb,G),
\end{equation}
where
$$\begin{cases}
    \boldsymbol{h}_{1}=\vb &\\
    \boldsymbol{h}_{j+1}=F_{j}(\boldsymbol{h}_{j},\eb,G) & j=1,\dots,m-1,
\end{cases}$$
so that the final output is obtained by applying the blocks $F_{1},\dots,F_{m}$ iteratively. We highlight how each message-passing-step transforms the node features but not the edge features. We also point out that, since a single message-passing-block allows neighbouring nodes to exchange information, a processor module with $m$ units allows communication between nodes that are $m$ edges faraway. Thus, by changing the number of message-passing-steps, one can move from local to nonlocal transforms (with the latter possibly being more expressive). However, we must also mention that large values of $m$ may give rise to oversmoothing phenomena\cite{rusch2023survey}, reason for which, in practice, a suitable compromise is required.

\subsubsection{Decoder module}
\label{sec:decoder}
In the end, the processor outputs some collection of node features $\vb':=\processor(\vb,\eb,G)$, with $\vb':V\to\mathbb{R}^{l}$. At this point, a terminal module, called the decoder, is exploited to recover the desired output. Here, we assume the latter to be consistent with the input signal $\ub$, and thus consist of $q$ nodal features. In practice, this is  achieved by relying on a suitable MLP unit $\Psi_{\decoder}:\mathbb{R}^{l}\to\mathbb{R}^{q},$ that transforms the original $l$ features onto the $q$ desired ones. In other words, the decoder module operates nodewise, and its action can be written as
\begin{equation}
\label{eq:decoder}
    \decoder(\vb'):=\Psi_{\decoder}\circ\vb'.
\end{equation}

\begin{figure*}
    \centering
    \includegraphics[width=0.95\textwidth]{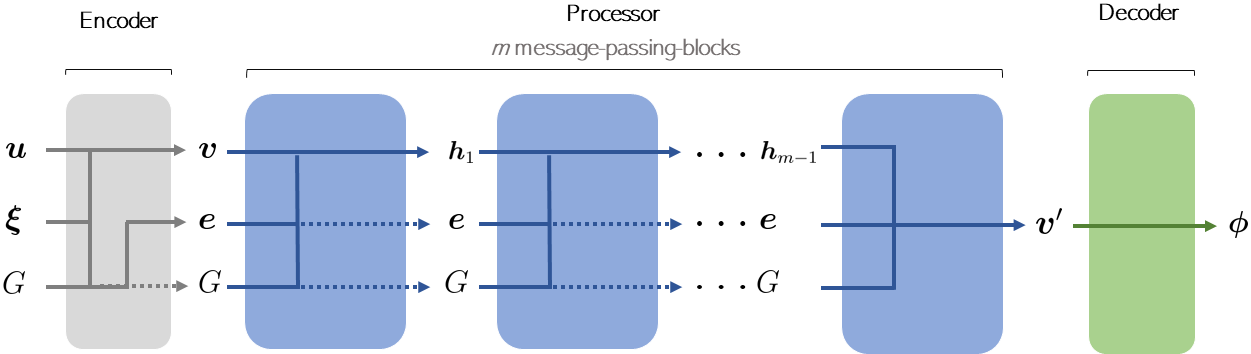}
    \caption{Visual representation of the Encoder-Processor-Decoder model, Section \ref{sec:epd}. Rigid arrows represent algorithmic computations (learnable and nonlearnable), while dashed arrows act as pointers (no computation implied). In gray, the encoder module, Eqs. \eqref{eq:encoder1}-\eqref{eq:encoder2}; in blue, the message-passing-blocks defining the processor unit, Eq. \eqref{eq:processor}; in green, the decoder module, Eq. \eqref{eq:decoder}.
    }
    \label{fig:epd}
\end{figure*}

\subsubsection{Overall architecture}

To summarize, the computational workflow of an Encoder-Processor-Decoder model reads 
\begin{equation}
\label{eq:epd}
\Phi(\ub,\xib,G):= \decoder(\processor( \encoder_{v}(\ub,\xib,G), \encoder_{e}(G), G)).\end{equation}
The reader can also find a visual depiction of Eq. \eqref{eq:epd} in Figure \ref{fig:epd}. Since the notation might be troublesome, we remark that the output of the Encoder-Processor-Decoder, $\boldsymbol{\phi}:=\Phi(\ub,\xib,G)$, is nothing but a collection of $q$-dimensional node features $\boldsymbol{\phi}:V\to\mathbb{R}^{q}.$

\begin{remark}
In the literature of surrogate and reduced order modeling, the words «Encoder» and «Decoder» are often associated to the concept of dimensionality reduction, where the two objects operate, respectively, to achieve data compression and reconstruction. Here, however, the meaning is completely different. The encoder module acts as a feature extractor and, in general, may increase the dimension of the input; conversely, the decoder module is used to map the nodal feature space onto the nodal output space (usually decreasing the dimension at each node), so to recover the quantities of interest. While the notation might be confusing to some of the readers, we have decided to stick to the one adopted by the GNN community \cite{pfaff2020learning,sanchez2020learning}.
\end{remark}

\section{\label{sec:surrogate}Application to surrogate modeling of parametrized PDEs}
Our goal is to predict an approximate solution $\Tilde{\mathbf{u}}$ at time $t^{n+1}$, given the state of the system at time $t^{n}$, for each node $i = 1, \ldots, \nfom$ of the computational mesh $\mesh$, that is
$$\Phi(\mathbf{u}^{n}_{\mub},t^{n},\mub)\approx\mathbf{u}^{n+1}_{\mub}.$$

Inspired by the general form of explicit Runge-Kutta methods, we model the time stepping scheme $\Phi$ by letting
\begin{equation}
\label{eq:time-stepping}
\Phi(\mathbf{v},t^{n},\mub):=\mathbf{v}+\Delta t\tilde{\Phi}(\mathbf{v},t^{n},\mesh),
\end{equation}
where $\tilde{\Phi}$ is a GNN architecture based on the Encoder-Processor-Decoder paradigm. Here, with little abuse of notation, we are identifying the computational mesh $\mesh$ with its underlying graph $G_{\mub}=(V_{\mub},E_{\mub})$, and the dof vector $\mathbf{u}^{n}_{\mu}\in\mathbb{R}^{\nfom}$ with its corresponding vertex feature map $\ub^{n}_{\mu}:V_{\mub}\to\mathbb{R}$. We also remark that in Eq. \eqref{eq:time-stepping}, the GNN model is made aware of the current time instant: in fact, according to our notation in Section \ref{sec:epd}, the latter is being interpreted as a global feature $\xi=t^{n}.$
%, that is learning the function $\Phi$ in (\ref{eq:epd}) such that
%\begin{equation}
%\left\{
%\begin{aligned}
%    &\Tilde{\ub}^{n+1}_i(\mub)  = \Phi(\Tilde{\ub}^{n}_{\mathcal{N}(i)},\xib_i,G ; \mub) \quad n = 0,\ldots,N_t - 1, \\
%    &\Tilde{u}^0_i(\mub)  = u_{h,i}^0(\mub),
%\end{aligned}
%\right.
%\label{eq:gnn}
%\end{equation}
%where $\mathcal{N}(i)$ is the set of the nodes in the neighborhood of node $i$. Equation (\ref{eq:gnn}) is the nodal version of Equation (\ref{eq:gnn_vec}), where the solution approximation at time $t^{n}$, i.e. $\Tilde{u}^{n}_{i}$, is contained in the \textcolor{red}{node features tensor $\mathbf{v}^n_{i}$} together with other relevant node features.\\
Therefore, the proposed approach aims at modeling the time-stepping scheme using an Encoder-Processor-Decoder architecture, $\tilde{\Phi}$, that incorporates both the nodal features at a specific time instance $t^n$, as well as the geometrical features of the mesh. The architecture aggregates information from neighboring nodes, processes it, and decodes the system solution at time $t^{n+1}$. This allows us to evaluate Equation \eqref{eq:time-stepping} independently of the number of dofs, while simultaneously accounting for the graph structure of the mesh. In particular, this makes it possible to train the model using a variety of different geometries and subsequently predict solutions for new meshes that were not included in the training data. 

This flexibility is guaranteed by the relational inductive bias of GNNs, which ultimately comes the message-passing paradigm: the model first computes the messages between neighboring nodes, and then performs a suitable aggregation of the information. %, in order to allow connections between nodes that are far from each other. Classical FFNNs are often prone to overfitting because they are fully connected, which means they may try to fit too closely the training data and fail to generalize to new examples. GNNs, on the other hand, can help to mitigate overfitting by incorporating information from neighboring nodes and edges, which can improve the robustness of the model. However, CNNs can also be used to limit overfitting, but as we previously discussed, the convolution operator may struggle with graph-based problems. Furthermore, 
In contrast, models based on FFNNs and CNNs are constrained by the number of nodes in the computational mesh, which prevents them from generalizing to different domains. The same issue is also encountered by other architectures, such as Mesh-Informed Neural Networks\cite{franco2022learning} (MINNs): in fact, even though MINNs can handle very complicated geometries at reduced training times, their implementation requires fixing the shape of the spatial domain and the resolution of the space discretization.
%Recently, a new class of architectures called Mesh-Informed Neural Networks (MINNs) has been introduced \cite{franco2022learning}. These architectures are specifically designed to handle mesh-based functional data and are capable of handling functional data defined on general domains of any shape better than classical FFNNs. Moreover, they can achieve reduced training times, lower computational costs, and better generalization capabilities. However, they require the number of nodes of the mesh $\nfom$ to be fixed, limiting their inductive capability.
In this sense, the additional flexibility provided by GNNs is extremely valuable.%For this reason, to increase flexibility in generalizing to different domains while varying the number of dofs, GNNs can be advantageous.

\subsection{Training and testing algorithms}
\label{sec:train}

From an operational point of view, the GNN model $\tilde{\Phi}$ in \eqref{eq:time-stepping} is trained on a suitable dataset of FOM solutions that serves as a ground truth reference. More precisely, after having constructed and initialized the GNN model, we exploit the FOM solver to generate a collection of training snapshots,
$$\{\mub_{i},\mathbf{u}_{\mub_{i}}^{0},\dots,\mathbf{u}_{\mub_{i}}^{N_{t}}\}_{i=1}^{\ntrain},$$
containing a total of $\ntrain$ different trajectories, each corresponding to a different geometrical configuration $\mub$. For the sake of simplicity, we assume that all the trajectories consist of $N_{t}$ snapshots in time: however, this assumption is not fundamental to our construction and it can easily dropped.

Let $\bm{\theta}$ be the vector collecting all the parameters of the GNN module. To emphasize the dependency of the latter on $\bm{\theta}$, let us write $\tilde{\Phi}_{\bm{\theta}}$ in place of $\tilde{\Phi}$. We train the GNN architecture by minimizing the
loss function below 
\begin{equation}
    \label{eq:loss}
    \begin{aligned}
    \mathcal{L}(\bm{\theta})  =  & cw_1\sum_{i=1}^{\ntrain}\sum^{N_{t}-1}_{n=0}|\mathbf{u}^{n+1}_{\mub} -  \mathbf{u}^{n}_{\mub}-\Delta t\tilde{\Phi}_{\bm{\theta}}(\mathbf{u}^{n}_{\mu_{i}},t^{n},\meshi)|^2\;+ \\
    &cw_2\sum_{i=1}^{\ntrain}\sum^{N_{t}-1}_{n=0}|\dot {\mathbf{u}}^{n}_{\mub}- \tilde{\Phi}_{\bm{\theta}}(\mathbf{u}^{n}_{\mu_{i}},t^{n},\meshi)|^2,
     \end{aligned}
\end{equation}
where $c=1/\ntrain N_{t}$ is a normalizing factor, whereas $w_{1}$ and $w_{2}$ are suitable hyperparameters to be tuned manually. %, and
% NOTA: in realtà, se si usa questo tipo di approssimazione, raccogliendo $1/\Delta t$ nel secondo termine della loss si ritrova esattamente il primo termine! Se vogliamo che abbia un senso, occorre essere un attimo più generali...
% Ora come ora, agisce solo come un ulteriore rescaling (probabilmente, Filippo ne stava traendo vantaggio perché, senza saperlo, questo rescaling lo aiutava a mitigare l'azzeramento dei gradienti).
%\begin{equation*}
%    \dot {\mathbf{u}}^{n}_{\mub} \approx \frac{\mathbf{u}^{n+1}_{\mub} - \mathbf{u}_{\mub} ^n}{\Delta t}
%\end{equation*}
%Propongo di metterla giù così
\nicola{The term $\dot{\mathbf{u}}^{n}_{\mub}$, instead, refers to}
%is
a \nicola{suitable} finite-difference approximation of the ground truth time-derivative \nicola{(e.g., computed by relying either on the forward or backward formulae)}. The loss function in \eqref{eq:loss} is made of two contributes: the first one, quantifies the error of the time-stepping scheme after a single iteration; the second one, instead, links the FOM derivative with the output of the GNN model. In particular, we do not rely on full rollouts or any other form of recursive training: this allows us to fully exploit the capabilities of GPUs tensor calculus and mitigate memory usage. 

Clearly, the downside to this is that, even after a successful training, our GNN model might be subject to error propagation when advancing in time multiple times. To limit this issue and ensure robust rollouts, we exploit the following strategies. At each epoch, that is, at each iteration of the optimization routine:
\begin{itemize}
    \item we do not directly optimize \eqref{eq:loss}, but rather rely on randomly selected mini-batches;
    \item we gitter the input data with random Gaussian noise, as to limit the sensitivity of $\tilde{\Phi}$ and to enhance the stability of the rollouts at prediction;
\end{itemize}

In practice, the optimization of the GNN model is carried out by relying on back-propagation \cite{rumelhart2986learning} and ADAM \cite{kingma2015adam}, with a variable learning rate that we decrease by a factor $\gamma>0$ after a specific number of epochs (see Algorithm \ref{alg:train}). In general, the training of the GNN model can be carried out iteratively until a stopping criterion is met. For instance, one may simply stop the training after a predefined number of epochs, see, e.g., Algorithm \ref{alg:train}.

\input{algorithm}

\noindent  
Once the model has been trained and a suitable vector of parameters $\bm{\theta}^{*}$ has been selected, the GNN is fully operational. That is, given any configuration of the geometric parameters $\mub$ and any initial condition $\mathbf{u}_{\mub}^{0}$, we can exploit the GNN model and \eqref{eq:time-stepping} \textit{online} to evolve the system iteratively and produce a complete rollout $\{\tilde{\mathbf{u}}_{\mub}^{n}\}_{n=0}^{N_{t}}$, where
\begin{equation}
\label{eq:gnn-dynamics}
\begin{cases}\tilde{\mathbf{u}}^{n+1}_{\mub}=\Phi(\tilde{\mathbf{u}}^{n}_{\mub},t^{n},\mub)& n\ge0\\
\tilde{\mathbf{u}}^{0}_{\mub}:=\mathbf{u}^{0}_{\mub}.&\end{cases}\end{equation}
Here, to further improve stability, one may also enforce any external constraint, such as Dirichlet boundary conditions, at each time iteration.

To test the quality of the GNN surrogate, we compare its predictions with those of the FOM for a set of new parameter instances. In particular, differently from the training stage, we now compare the overall trajectories and use the GNN to produce full rollouts of the solution.

Quantitatively, we compute the prediction error as the relative MSE (RMSE) error between the network prediction and the ground truth solution:
\begin{equation}
RMSE(\Tilde{\mathbf{u}}_{\mub}^{1},\dots,\Tilde{\mathbf{u}}_{\mub}^{N_{t}};\;\mathbf{u}_{\mub}^{1},\dots,\mathbf{u}_{\mub}^{N_{t}}) = \frac{1}{N_{t}}\sum_{n = 1}^{N_{t}}\frac{|\Tilde{\mathbf{u}}^n_{\mub}-\mathbf{u}^n_{\mub}|^2}{|\mathbf{u}^n_{\mub}|^2},
\end{equation}
where the GNN rollout is obtained as in \eqref{eq:gnn-dynamics}.

\section{\label{sec:exp}Numerical experiments}
\input{experiments}

\section{\label{sec:conclusions}Conclusions}
\nicola{We presented a novel approach to surrogate modeling based on Graph Neural Networks (GNNs) for the efficient evolution of dynamical systems defined over parameter dependent spatial domains. The approach differs substantially from classical Reduced Order Modelling techniques, in that it provides a way to handle parameter dependent PDEs with a variable number of degrees of freedom.
The method is based on a data-driven time-stepping scheme that explicitly accounts for the Markovian structure of the dynamical system, while also including geometric information via GNN modules. The approach is shown capable of yielding stable simulations, even for long rollouts, while simultaneously generalizing to unseen geometries, thus providing remarkable benefits when compared to other techniques based on different neural network architectures.}

\nicola{Although limited to fairly simple problems, our results indicate that GNNs can be a valuable tool for ROM practitioners, providing researchers with new ways for handling geometric variability.
Future research may involve the exploration of hybrid approaches where GNNs are combined with other well-established Deep Learning-based reduced order models, such as autoencoders and U-Net-like architectures, in an attempt to generalize the whole idea %his could be particularly interesting as it would allow one to extend the approach %This hybrid approach could 
to more complicated problems with thounsands or millions of degrees of freedom.}

%potentially lead to improved accuracy in performing simulations, even on more refined meshes, while simultaneously reducing computational complexity.

%Furthermore, additional investigations
Another interesting  question could be whether this approach can benefit from the integration of suitable %explore the use of 
attention mechanisms, or other forms of neural network architectures, that can selectively weight the contributions of different nodes in the graph. We leave these considerations for future work. %, potentially improving the performance of the model on more complex geometries and physical phenomena.

\section*{Acknowledgments}
The present research is part of the activities of project Dipartimento di Eccellenza 2023-2027, funded by MUR, and of project FAIR (Future Artificial Intelligence Research) project, funded by the NextGenerationEU program within the PNRR-PE-AI scheme (M4C2, Investment 1.3, Line on Artificial Intelligence). NF, SF and AM are members of Gruppo Nazionale per il Calcolo Scientifico (GNCS) and of Istituto Nazionale di Alta Matematica (INdAM).

\section*{Data Availability Statement}
The data that support the findings of this study are available from the corresponding author upon reasonable request.

\section*{Author Declarations}
The authors have no conflicts to disclose.

%\appendix

%\section{Appendixes}

\bibliographystyle{plain}
\bibliography{aipsamp}% Produces the bibliography via BibTeX.

\end{document}

%% file: algorithm.tex
\begin{algorithm}[H]
\vspace{0.5em}\textbf{Input:} network $\tilde{\Phi}$; timestep $\Delta t$; a list of $N_{train}$ training trajectories $\mathbf{U}$ (each of length $N_{t}$); a list of edge connectivity matrices $\mathbf{E}$; a list of edge features matrices $\mathbf{W}$; a list of inner nodes $\mathbf{I}$; learning rate $\nu$; decay factor $\gamma$; number of training epochs $\texttt{epochs}$; batch size $N_{b}$; noise variance $\sigma^2$.\vspace{0.5em}\\
\textbf{Output:} optimal model parameters $\bm{\theta}^*$.\vspace{0.5em}
\caption{Training Algorithm}
\label{alg:train}
\begin{algorithmic}[1]
\STATE $\texttt{epoch} = 0$.
\STATE Randomly initialize $\bm{\theta}^0$.
\WHILE{$\texttt{epoch} < max\_epoch$} 
    \STATE Create the list $indices = [1,\ldots,N_{train}]$ and shuffle it randomly.
   \FOR {$sim$ in $indices$} 
    \STATE    $\mathbf{U}_{sim} = \mathbf{U}[sim]$, $\mathbf{U}_{sim} \in \mathbb{R}^{N_t \times N_h \times q}$ where $N_t$ is the total number of time instances, $N_h = N_h(\bm{\mu}_{sim})$ are the mesh dofs and $q$ is the number of node features.\vspace{0.5em}
    \STATE    $\mathbf{E}_{sim} = \mathbf{E}[sim]$, $\mathbf{E}_{sim} \in \mathbb{R}^{N_{edges} \times 2}$.\vspace{0.5em}
    \STATE    $\mathbf{W}_{sim} = \mathbf{W}[sim]$, $\mathbf{W}_{sim} \in \mathbb{R}^{N_{edges} \times N_e}$ where $N_e$ is the number of edge features.\vspace{0.5em}
    \STATE $\mathbf{I}_{sim} = \mathbf{I}[sim]$, $\mathbf{I}_{sim} \in \mathbb{R}^{N_h}$ with $\mathbf{I}_{sim}[i] = 1$ if node $i$ is an inner node, $0$ otherwise.\vspace{0.5em}
    \STATE    $b = 0$.\vspace{0.5em}
    \WHILE{ $b < N_t$ }
        \STATE    $\mathbf{U}_{b} = U_{sim}[b:b + N_{b}]$.
        \STATE    Create noise tensor $\bm{\Sigma} =  \sigma \mathbf{Z}$ where $\mathbf{Z} \in \mathbb{R}^{N_{b} \times N_i \times q}$ is a random tensor with $N_i$ inner nodes.
        \STATE Initialize $\mathbf{U}_{noise} = \mathbf{U}_{b}$.
        \STATE    $\mathbf{U}_{noise}[:,I_{sim}] \ += \ \Sigma$.
        \STATE     Calculate target derivative $\mathbf{U}_{dot} = (\mathbf{U}_{b}[1:] - \mathbf{U}_{noise}[:-1])/\Delta t$.
        \STATE    Make a forward pass through the network $\tilde{\Phi}(\mathbf{U}_{noise}, \mathbf{E}_{sim}, \mathbf{W}_{sim})$.
        \STATE    Calculate network solution $U_{net} = \mathbf{U}_{noise}[:-1] + \Delta t \tilde{\Phi} $.
        \STATE    Calculate training loss $\mathcal{L}_{b}$.
        \STATE    Back-propagation through the net and parameters update: $\bm{\theta}^1 = \textnormal{ADAM}(\nu,\bm{\theta}^0$).
        \STATE $\bm{\theta}^0 = \bm{\theta}^1$.
        \STATE    $b \gets b + N_{b}$.
    \ENDWHILE
    \ENDFOR\vspace{0.5em}
    \IF{$\text{mod}(\texttt{epoch}, 500)=0$}
    \STATE Reduce learning rate by a factor $\gamma$.
    \ENDIF
    \STATE $\texttt{epoch} \gets \texttt{epoch} + 1$.
\ENDWHILE\vspace{0.5em}\\
Pick the last weights updated $\bm{\theta}^1$.\vspace{0.5em}
\end{algorithmic}
\end{algorithm}

%% file: experiments.tex
In this Section, \nicola{we assess the capabilities of the proposed approach over three advection-diffusion problems} %shall now test the proposed approach } the technique introduced in the previous sections is applied to an 
%advection-diffusion problem in three different cases, 
of increasing complexity: 
\begin{itemize}
    \item \nicola{a scalar diffusion in} a 2D square with a circular obstacle \nicola{and} a time-varying advection term;
    \item a 2D Stokes flow in proximity of a bump;
    \item a 3D Stokes flow around a cylinder.
\end{itemize}

All the examples \nicola{are characterized by parameter dependent spatial domains, where a given obstacle is allowed to move across the domain, with possible changes in terms of shape and dimension.} %are studied in domains in which we have let the obstacle vary its position; then, the first two examples are also studied in domains where the obstacle varies both its position and dimension. In this way, we have exploited the model capacity to infer the geometrical features of the mesh.
\nicola{In this way, we can effectively test the ability of GNNs in handling geometric variability.}

\subsection{Advection-Diffusion problem in a square domain with a circular obstacle}
%\subsubsection{Problem definition}
\label{sec:ad}
%In this section, the model is applied to the 
\nicola{To start, we consider the following} advection-diffusion problem:
\begin{equation}
\label{eq:ad}
\begin{cases}
    \displaystyle
    \frac{\partial u}{\partial t} - D\Delta u + \mathbf{b}\cdot\nabla u= 0&\text{in} \ \Omega \times (0,T] \\
    u(x,y) = (x-1)^2+(y-1)^2 & \text{on} \ \partial\Omega \times (0,T]\\
    u_{0}(x,y) = (x-1)^2+(y-1)^2 & \text{in} \ \Omega,
\end{cases}
\end{equation}
where $\Omega = (0,1)^2 \setminus C$, with $$C = \{(x,y) : \ (x-c_{x})^2 + (y-c_{y})^2 \leq (0.15)^2\}.$$ 
Here, we set $T = 2$, $D = 0.1$ and $\mathbf{b}(t)=[1-t,1-t]$. \nicola{In particular, due to the time-varying convection field, $\mathbf{b}$, the resulting dynamical system can be regarded as nonautonomous.} % allows us to observe a more complex phenomenon during the time horizon considered with respect to the case in which $\mathbf{b}$ is a steady vector field. 
In our simulations, we parametrize the center of the circle as $$\mathbf{\mub} = (c_{x},c_{y}) \in \Theta := \{(x,y) : \ 0 < x < 1, \ y \geq 0.5\},$$ \nicola{which we let vary as we generate the training data.}
In agreement with Equations \eqref{eq:pde}-\eqref{eq:fem}, the ground truth FOM simulations of Problem (\ref{eq:ad}) are obtained by first discretizing in space \nicola{via P1 Continuous Galerkin Finite Elements, and then in time using the Backward Euler Method.} %Space and imposing Dirichlet boundary conditions. Then, we discretize in time using the Backward Euler Method following Equation (\ref{eq:fem}). %Two examples of FOM solutions are reported in Figure \ref{fig:ad}.\\
%\begin{figure*}
%    \centering
%  \includegraphics[scale=0.6]{Images/ad_1.png}
%  \includegraphics[scale=0.6]{Images/ad_2.png}
%  \caption{Test case 1, Advection-Diffusion problem, FOM solutions. Both rows represent 3 time steps of simulation. In the first row, the center of the obstacle is at $\mub = (0.5,0.5)$, while in the second row is at  $\mub = (0.32,0.68)$. }
%  \label{fig:ad}
%\end{figure*} 
\nicola{The time step chosen is $\Delta t = 0.02 $, resulting in 101 time snapshots for each simulation. We also mention that, following our notation in Section \ref{sec:gnns}, we have $q=1$, as the solutions to \eqref{eq:ad} are scalar fields.}

\subsubsection{Problem data}
We collected a dataset composed of 100 simulations, each obtained for a different position of the center of the obstacle, with a number of mesh nodes varying from 770 to 790. The training set is composed of 80 \nicola{randomly selected}  simulations, while the \nicola{remaining 20 are kept for testing}. %test set has 20 simulations, both chosen randomly among the 100 FOM simulations run. 

\nicola{For what concerns the design of the GNN architecture and its training, we have reported a synthetic overview in Table \ref{tab:hyperparameters}. In particular, in this case, we adopt a simplified loss function that only features the approximation of the time-derivative; in other words, we set $w_1 = 0$ and $w_2 = 1$ in Equation \eqref{eq:loss}.}

%For each batch, the loss is computed using only the MSE between the derivative computed by the network and the ground truth derivative calculated with the finite differences scheme presented in Section \ref{sec:train}. In this case, we set the loss weights as $w_1 = 0$ and $w_2 = 1$, so we do not have the term that involves Forward Euler prediction in the loss function.

\renewcommand{\arraystretch}{1.5}
\begin{table*}
    \centering
    \begin{tabular}{|l|llll|lllllll|}
    \hline
    \textbf{Problem}  &  \textbf{MP} & \textbf{MLP} &  $l$ &  \textbf{Activ}&  \textbf{Epochs}&  \textbf{Batch} &\textbf{lr}&  $\gamma$ & $\sigma^{2}$ & $w_{1}$ & $w_{2}$\vspace{-0.5em}\\
    &\textbf{steps}&\textbf{layers}&&&\textbf{(max)}&\textbf{size}&&&&&\\
    \hline
    %\hline
    \textbf{Adve-diff} & 12 & 2 & 32 & SiLU & 1500 & 25 & $10^{-3}$ & 0.1 & $10^{-6}$ & 0 & 1\\   
    \textbf{Stokes 2D} & 18 & 2 & 32 & SiLU & 3000 & 25 & $10^{-3}$ & 0.1 & $10^{-6}$ & 0.5 & 0.5\\  
    \textbf{Stokes 3D} & 15 & 2 & 32 & SiLU & 2000 & 25 & $10^{-4}$ & 0.1 & $10^{-5}$ & 0.5 & 0.5\\  
    \hline
    \end{tabular}
    \caption{GNN architecture and training hyperparameters for the three case studies. MP = message passing, MLP layers = (common) depth of all the MLP units in the Encoder-Processor-Decoder pipeline, $l$ = (local) feature space dimension, lr = learning rate, $\gamma$ = learning rate decay factor (applied every 500 epochs), $\sigma^{2}$ = noise variance, $w_{i}$ loss function weights. SiLU = Sigmoid weighted Linear Unit, $x\to x/(1+\exp(-x)).$}
    \label{tab:hyperparameters}
\end{table*}

\subsubsection{Numerical results}
%The results of the rollout predictions of the test simulations are summarized 
\nicola{Results are} in Table \ref{table:errors}. As we can see, all the predictions RMSEs are of order $10^{-3} $ to $10^{-4}$. Moreover, our model outperforms significantly the ground truth solver in the simulation time at testing stage. %The best prediction, together with its ground truth solution, is shown in Figure \ref{fig:best}. \\
The dynamic of the problem is well predicted and no propagation errors are spotted. Hence, our model %and the training strategy adopted seem to be in principle good at 
\nicola{appears capable of} solving problems concerning evolutionary PDEs, \nicola{in that it can approximate multiple time steps in a stable way.} %in which multiple %rollout
%time steps need to be predicted.

\begin{figure*}
\centering
  \includegraphics[scale=0.55]{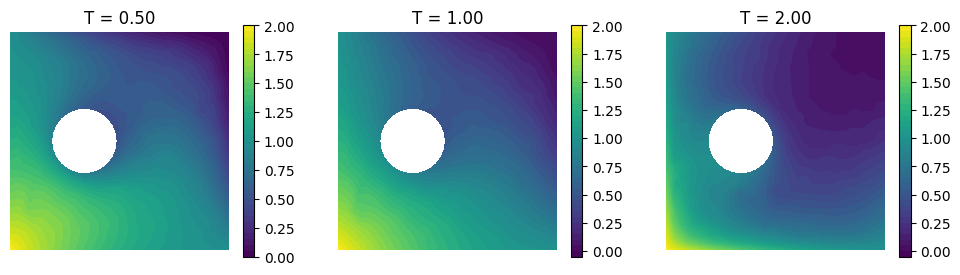}
  \includegraphics[scale=0.55]{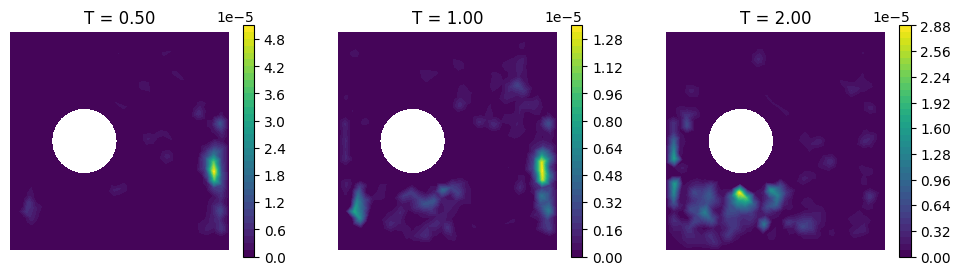}
  \caption{\textnormal{Test case 1, Advection-Diffusion problem. Prediction obtained for $\mub = (0.29,0.5)$ with the obstacle close to the source. First row: rollout prediction. Second row: RMSE related to each time step between the prediction and the corresponding ground truth solution.} }
  \label{fig:pred5}
\end{figure*}

\begin{figure*}
 \centering
  \includegraphics[scale=0.55]{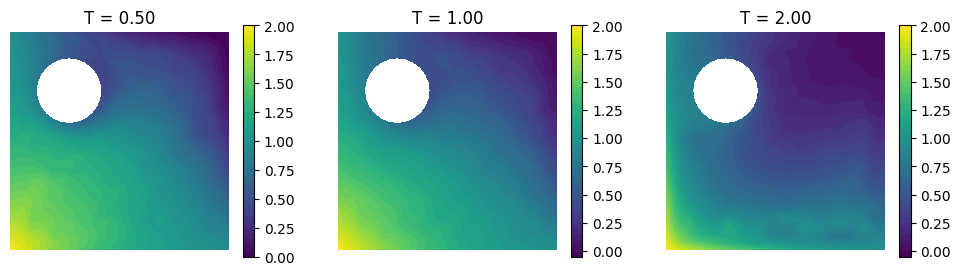}
  \includegraphics[scale=0.55]{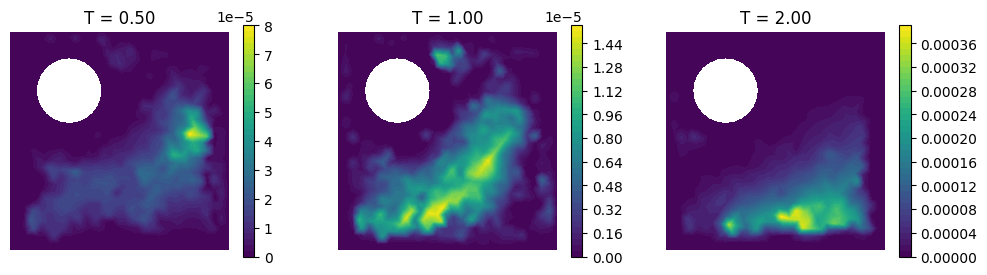}
  \caption{\textnormal{Test case 1, Advection-Diffusion problem. Prediction obtained for $\mub = (0.25,0.75)$ with the obstacle on the top left corner. First row: rollout prediction. Second row: RMSE related to each time step between the prediction and the corresponding ground truth solution.} }
  \label{fig:pred19}
\end{figure*}

\nicola{Still, it is worth looking at some of the simulations obtained during the testing phase, as to further appreciate the ability of the proposed approach in handling different geometric configurations. For instance, Figures \ref{fig:pred5} and \ref{fig:pred19} show two different GNN rollouts corresponding to two different positions of the obstacle. Despite these trajectories being different from the ones seen during training, the model manages to capture all the main features characterizing the solutions, such as the behavior near the obstacle and the direction of propagation.}

We highlight that a GNN-based approach follows a \textit{local-to-global} \nicola{paradigm,} %generalization approach, 
\nicola{first processing information at the node level (encoder), and then aggregating the output at the neighbour level (processor). Clearly, the lack of smoothness in PDE solutions can pose some challenges, as GNNs are known to struggle with capturing such properties. In this sense, it is not surprising to see that the prediction in Figure \ref{fig:pred19} is worse than the one in Figure \ref{fig:pred5}. In fact, in the former case, the obstacle is closer to the corner of the spatial domain. Of note, we mention that the trajectory in Figure \ref{fig:pred19} is actually the worst across the whole test set.} 

\nicola{It is also interesting to see that the prediction error exhibits an oscillating trend. In fact, after a first increase, the accuracy appears to improve ($t=0.5$ vs $t=1.00$), which is most likely caused by the presence of a diffusion phenomenon; then, however, the approximation deteriorates again due to the presence of the convection field, which pushes the errors either towards the obstacle (Figure \ref{fig:pred5}) or the bottom boundary (Figure \ref{fig:pred19}).}

% what has been processed \nicola{within the neighborhood}. from the neighborhood to connect the mesh nodes. However, this approach may result difficult for some trajectories and yield worse predictions at critical points of the mesh, such as along the boundary. The lack of smoothness in PDE solutions is a common challenge faced by deep learning-based surrogate models, as neural networks are known to struggle with capturing such properties. \par

\begin{figure}%
    \centering
    \includegraphics[width=0.495\textwidth]{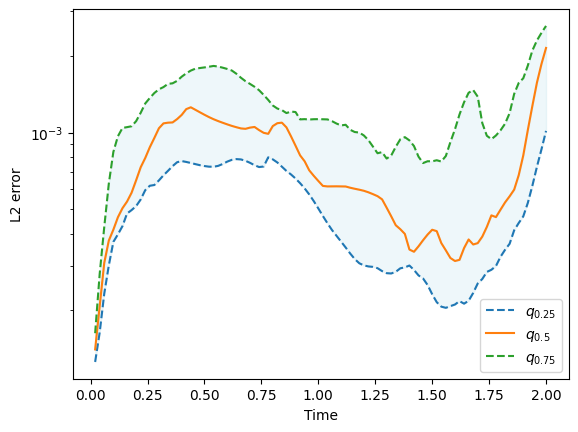} %
    \includegraphics[width=0.495\textwidth]{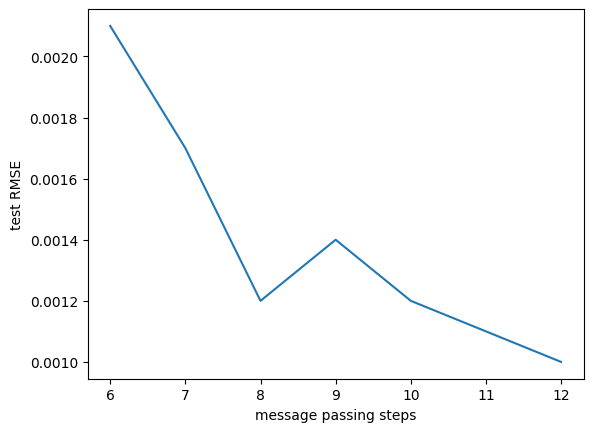} %
    \caption{\textnormal{Test case 1, Advection-Diffusion problem. Left: $L^2$ relative error vs Time plot. The dashed lines represent the first and the third quantiles of the $L^2$ errors among all the test predictions, while the orange line is the median. The shaded area can be considered a confidence region for the simulation error. Right: Test RMSE vs message passing steps.}}%
    \label{fig:mpsteps}%
\end{figure}

\nicola{We can further appreciate this phenomenon in Figure \ref{fig:mpsteps} (top row), where we have synthesized the dynamics of the relative $L^{2}$-error. More precisely, the picture shows how the quality of the approximation changes within time: to account for the variability in the test set, both median and quantile curves are reported. The trend appears to be fairly general, although the behavior quickly differentiates among different simulations (note how, as $t$ goes from 0 to 0.25, the upper and lower quartiles rapidly split apart).
} %Another key aspect to analyze is the trend of the $L^2$ relative error as the simulation time varies. This error represents, at each time step, the $L^2$ error between the prediction and the FOM solution divided by the $L^2$ norm of the FOM solution. In Figure \ref{fig:mpsteps}, on the left, we can see 3 curves that represent the evolution of the $L^2$ error between our solution and the ground truth over the time interval $[0,T]$. The curves show respectively how the first, second, and third quantiles of the error vary in time. The trend is coherent with the results previously shown, indeed we can notice a substantial increase in the error as the simulation time increases, and in particular, the final instants strongly influence the RMSE of the predictions. The first quantile and the median are similar, while the third shows more oscillations in the error. However, for most of the simulation time, we have good bounds for the error, which is kept around $10^{-3}$. %\par
%Despite the anomaly shown in Figure \ref{fig:pred19}, it is worth noticing that the model is completely mesh-independent, and for the majority of trajectories provides a good prediction for long rollouts, keeping the number of parameters to a surprisingly small number. Indeed for this example, the model has 91009 parameters, which is easy to reach if you try to solve a similar problem with a fully connected network.

\subsubsection{The \textit{message passing steps} hyperparameter}

Among all the hyperparameters, the one most influencing the goodness of the model is the number of message-passing steps. This number represents how much in-depth we look at the neighborhood when we propagate the message. A small number of message-passing steps may result in underfitted areas of the mesh, while a big one will slow down the training, increasing too much the number of parameters, possibly yielding overfitting.\par

\nicola{Here, we tuned this parameter via trial and error.} A plot of the corresponding results %found 
can be seen in Figure \ref{fig:mpsteps} (bottom row). The test RMSE reaches a local minimum for $m=8$ message-passing steps. This is the best choice if we want to keep control of the number of total parameters of the network, which are only 61825 in this case.
\nicola{However, since the architecture obtained for $m=12$ is still reasonably complex, we stick to the latter one. We do not proceed further as the improvement rate, in terms of $m$, no longer justifies favoring a larger number of message-passing steps.} %, we the overall performance improves for a higher number of message-passing steps. Here, we stop as $m=12$, as the keep we have a better approximation of the trajectories for almost all the simulations and the number of parameters does not grow too much.

\renewcommand{\arraystretch}{1.5}
\begin{table}
    \centering 
    \begin{tabular}{| l l l l l l|}
    \hline
     & \textbf{RMSE} \;\;& \textbf{RMSE} & \textbf{RMSE} \;\;& $t_{\text{FOM}}$ & $t_{\text{GNN}}$\vspace{-0.5em}\\
     &(mean) & (max) & (min) & &\\
    \hline 
    \;\textbf{Adve-diff} §\ref{sec:ad}\;\; & $1.20\mathrm{e}{-3}$ & $6.10\mathrm{e}{-3}$ & $4.0\mathrm{e}{-4}$ &$159.80s$  & $9.83s$\;\;\\
    \;\textbf{Stokes 2D} §\ref{sec:stokes}\;\; & $1.64 \mathrm{e}{-2}$ & $7.35 \mathrm{e}{-2}$ & $1.2 \mathrm{e}{-3}$ &$115.65s$  & $ 7.51s$\;\;\\
    \;\textbf{Stokes 3D} §\ref{sec:stokes3d}\;\; & $4.37 \mathrm{e}{-2}$ & $6.24 \mathrm{e}{-2}$ & $1.9 \mathrm{e}{-2}$ &$729.42s$\;  & $10.4s$\;\;\\
    \hline
    \end{tabular}
    \\[10pt]
    \caption{Comparison between FOM and GNN-surrogate in terms of model accuracy and computational time for the three case studies.}
    \label{table:errors}
\end{table}

\subsubsection{Generalization to obstacles with different dimensions}
%This advection-diffusion example 
\nicola{Problem \ref{eq:ad}} can also be extended to domains in which both the position and the dimension of the obstacle change. \nicola{To this end}, we modify our training dataset slightly by adding new simulations in which the obstacle has either a smaller and a larger radius% than before and this will improve the generalization of the model. 
\nicola{Mathematically speaking, this corresponds to considering an augmented parameter space where} % the geometrical parameter now is 
$\mub = (c_x,c_y,r) \in \Theta = \{(x,y): \ 0 < x < 1, \ y \geq 0.5\} \times \{0.1, 0.15, 0.2 \}$. 

\begin{figure*}
 \centering
  \includegraphics[scale=0.55]{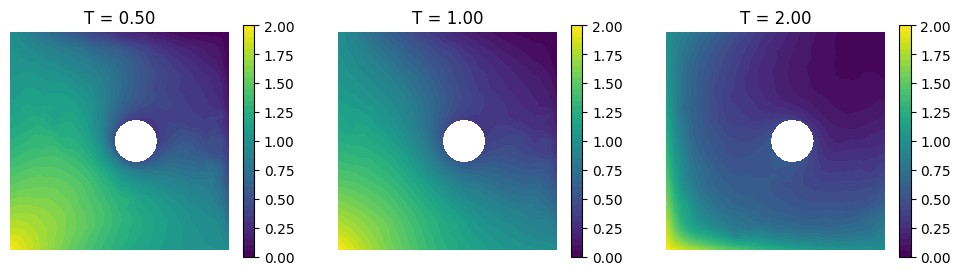}
  \includegraphics[scale=0.55]{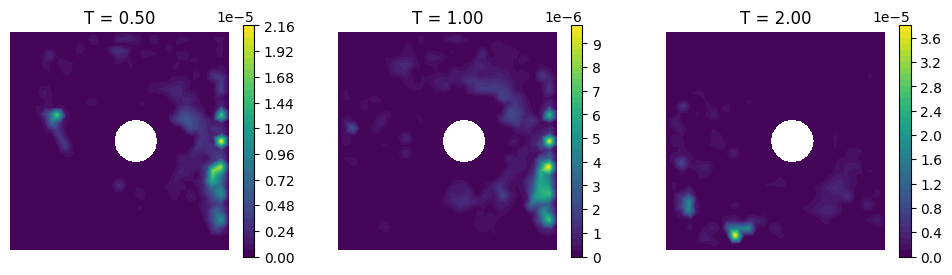}
  \caption{\textnormal{Test case 1, Advection-Diffusion problem. Prediction obtained for $\mub = (0.6,0.52,0.1)$. First row: rollout prediction. Second row: RMSE related to each timestep between the prediction and the corresponding ground truth solution.}}
  \label{fig:pred13obs}
\end{figure*}

\begin{figure*}
 \centering
  \includegraphics[scale=0.55]{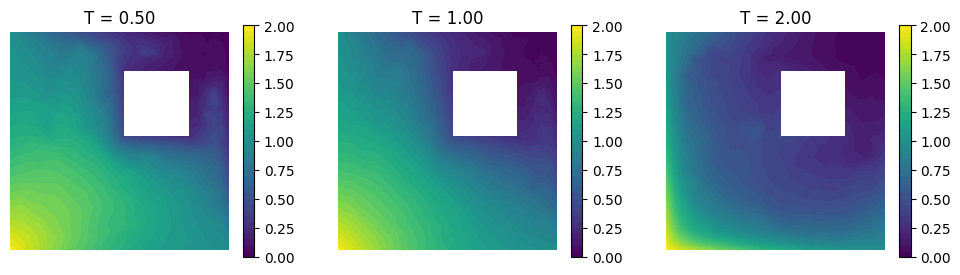}
  \includegraphics[scale=0.55]{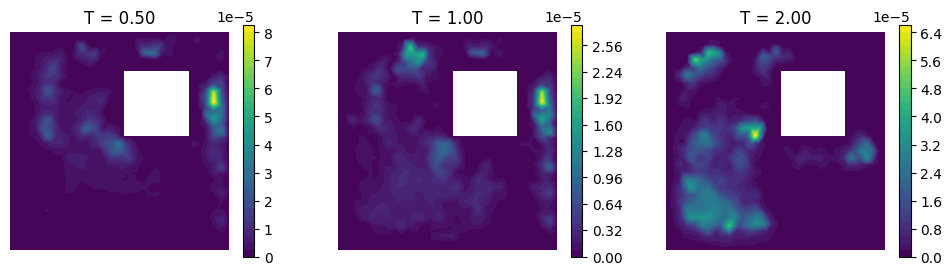}
  \caption{\textnormal{Test case 1, Advection-Diffusion problem. Prediction obtained for $\mub = (0.7,0.7,0.3)$ with a square obstacle. First row: rollout prediction. Second row: RMSE related to each timestep between the prediction and the corresponding ground truth solution.}}
  \label{fig:predsq}
\end{figure*}

Again, we test the model on new simulations which have varying obstacle positions and dimensions. In Figure %\ref{fig:pred7obs} and 
\ref{fig:pred13obs} the prediction for a new test simulation is reported. The model can generalize well on this problem even if the geometries differ a lot from each other, in terms of sizing. Moreover, there is no need to increase the number of message-passing steps, \nicola{meaning that the GNN architecture has the same complexity as before.} %so the number of parameters can be kept under control.
%\begin{figure*}
% \centering
%  \includegraphics[scale=0.55]{Images/pred_7_obstacle.png}
%  \includegraphics[scale=0.55]{Images/gt_7_obstacle.png}
%  \caption{\textnormal{Test case 1, Advection-Diffusion problem. Prediction obtained for $\mub = (0.55,0.53,0.2)$ . First row: rollout prediction. Second row: ground truth solution.}}
%  \label{fig:pred7obs}
%\end{figure*} 

\;\\An important question that arises is whether our model can predict solutions where the obstacle has a different shape, \nicola{and whether we can achieve this without having to retrain the whole network.} %without requiring retraining of the network. 
To investigate this, we present an example in Figure \ref{fig:predsq} of a prediction obtained with a square obstacle located in the top right of the domain\nicola{: strictly speaking, this configuration cannot be described in terms of our previous parametrization; nonetheless, we can still apply our GNN surrogate, as the latter only depends on the geometrical parameters through the underlying mesh (by itself, the parametrization never enters the equation).} % that is corresponding to the parameters $\mub = (c_x,c_y,L) = (0.7,0.7,0.3)$, where $r_x$ and $r_y$ are the coordinates of the center of the square, and $L$ is the length of its edge.

Surprisingly, the errors are of the same order of magnitude as those discussed earlier, and the prediction of the overall dynamics is remarkably accurate. This result is attributed to the ability of the model to understand different geometries by means of its inductive structure. GNNs, in particular, can automatically incorporate the geometrical structure of the domain by utilizing both the edge connectivity matrix and the edge features. However, some difficulty is observed in handling the nodes surrounding the obstacle, especially at the corners, but this does not appear to affect the overall accuracy of the prediction. These findings suggest that our model has the potential to generalize well to other geometries, without the need for extensive retraining, thus enhancing its practical applicability in real-world scenarios.

\subsection{Advection-Diffusion problem in a 2D Stokes flow in proximity of a bump}
\label{sec:stokes}
We now consider another advection-diffusion problem as \eqref{eq:ad}, where the advection field $\mathbf{b}$ is no longer fixed by hand, but it is rather obtained by solving the following stationary Stokes problem:
\begin{equation}
\label{eq:stokes}
\begin{cases}
-\nu \Delta \mathbf{b} + \nabla p &= 0 \qquad \text{in} \ \Omega\\
\nabla \cdot  \mathbf{b} &= 0 \qquad  \text{in} \ \Omega
\end{cases}
\end{equation}
where $p$ is the pressure field and the boundary conditions are given by:
%\begin{align}
%\mathbf{b}  &= 0 \qquad \text{on} \ \Gamma_{D} \\
%\mathbf{b} &= \mathbf{b_{in}} \quad  \text{on} \ \Gamma_{in}\\
%\nu \frac{\partial \mathbf{b}}{\partial \mathbf{n}} - p\mathbf{n} &= 0 \qquad \text{on} \ \Gamma_{N},
%\end{align}
\begin{equation*}
\mathbf{b}  = 0 \ \text{on} \ \Gamma_{D},\quad
\mathbf{b} = \mathbf{b_{in}} \ \text{on} \ \Gamma_{in},\quad
\nu \frac{\partial \mathbf{b}}{\partial \mathbf{n}} - p\mathbf{n} = 0 \ \text{on} \ \Gamma_{N},
\end{equation*}
with
\begin{align}
\mathbf{b_{in}} = \left(\frac{40Uy(0.5-y)}{0.5^2},0\right), \quad U = 0.3, \quad \nu = 10^{-3};
\end{align}
$\mathbf{b_{in}}$ represents the value of $\mathbf{b}$ at the inflow $\Gamma_{in}$, while $\Gamma_{D}$ and $\Gamma_{N}$ denote the Dirichlet wall (top and bottom) sides, and the Neumann right outflow boundaries, respectively.  

\nicola{This time,} the domain $\Omega$ is a rectangular channel $(0,1) \times (0,0.5)$ with a parametrized bump along the top wall edge. Here $\Gamma_{in} = \{ x = 0 \}$, $\Gamma_{D} = \{ y = 0 \} \cup \{ y = 0.5 \}$ and $\Gamma_{N} = \{ x = 1 \}$. During our simulations, we shift the position of the bump in a way that its center $c_x$ varies from $0.35$ to $0.65$. Hence, we consider $\mu \in \Theta = [0.35,0.65]$.

Regarding the advection-diffusion problem, at the inflow $\Gamma_{in}$ we impose the Dirichlet boundary condition $$u_{in}(x,y) = \frac{4y(0.5 - y)}{0.5^2},$$ which is also the initial condition, while on $\Gamma_{D}$ we impose no-slip boundary conditions and on $\Gamma_{N}$ we set $\partial u/\partial n = 0$. The final simulation time is $T = 0.5$ and $D = 0.01$. Our results will only focus on the approximation of the solution $u$ of the advection-diffusion problem, despite the latter also depending implicitely on \eqref{eq:stokes}. %depends on the parameter $\mu$ implicitly through the solution of the Stokes problem.  % and an example of a FOM solution is reported in Figure \ref{fig:bumpgt0}. \\
%\begin{figure*}
%\centering
%  \includegraphics[scale=0.7]{Images/bump_gts_0.png}
%  \caption{\textnormal{Test case 2, Advection-Diffusion problem in a 2D Stokes flow. 3 timesteps of a FOM simulation. The center of the bump here is at $\mu = 0.4696$.} }
%  \label{fig:bumpgt0}
%\end{figure*} 

\subsubsection{Problem data}
Our dataset is composed of 125 simulations, each obtained for a different position of the bump. In each of these cases, the mesh is rebuilt yielding a number of mesh nodes varying from 937 to 1042. The chosen time step is $\Delta t = 0.01$, resulting in 51 time snapshots for each simulation. %We consider the same node and edge input features of the previous example. However, compared to the previous test case, the loss function is slightly different, 
%since now the two terms defined in Section \ref{sec:train} are weighted equally 
The training set is made by 100 simulations, while the test set includes 25 simulations, both chosen randomly among the 125 FOM simulations.

\nicola{Differently from our previous test case, we consider a loss function where the two terms in \ref{eq:loss} are weighted equally (cf. Table \ref{tab:hyperparameters}). As before, we refer to Table \ref{tab:hyperparameters} for further details about GNN and training hyperparameters.}
%Most of the hyperparameters are chosen equally to the previous test case of Section \ref{sec:ad}, except for the following ones: %that is: 

%\begin{itemize}
%\item $l = 32$, where $l$ is the latent dimension of the network;

%\item noise variance $\sigma^2 = 10^{-6}$;

%\item total number of epochs $max\_epoch = 3000$;

%\item learning rate $\nu = 10^{-3}$ with decay $= 0.1$ after 500 and 1000 epochs;

%\item SiLU is chosen as activation function for each MLP layer;

%\item number of MLP layers $mlp\_layers = 2$;

%\item message passing steps $mp\_steps = 15$.
%\end{itemize}

\subsubsection{Numerical results}
As before, quantitative results are in Table \ref{table:errors}.
%The results of the rollout predictions of the test simulations are summarized in Table \ref{table:errors}.
In this more complex problem, the RMSEs are higher than the ones obtained in the previous example; however, \nicola{the predictions are still fairly accurate and} we still outperform the ground truth solver in terms of time efficiency. %\par

\begin{figure*}
\centering
  \includegraphics[scale=0.6]{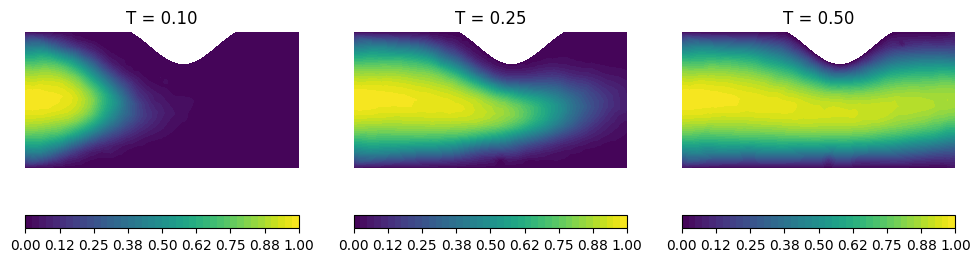}
  \includegraphics[scale=0.6]{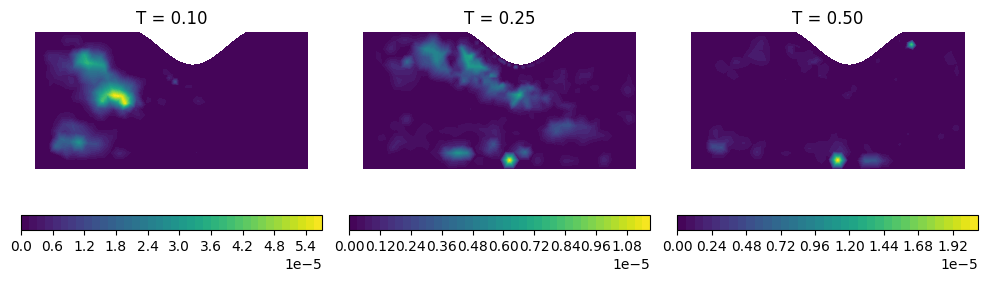}
  \caption{\textnormal{Test case 2, Advection-Diffusion problem in a 2D Stokes flow. Prediction obtained for $\mu = 0.58$ with the bump on the right part of the upper edge. First row: rollout prediction. Second row: RMSE related to each timestep between the prediction and the corresponding ground truth solution}}
  \label{fig:pred12bump}
\end{figure*}

\begin{figure*}
\centering
  \includegraphics[scale=0.6]{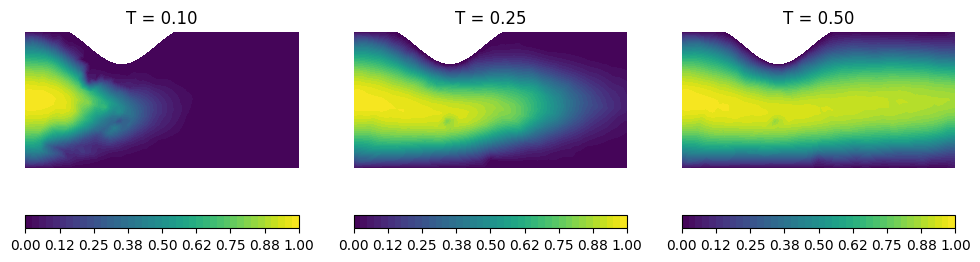}
  \includegraphics[scale=0.6]{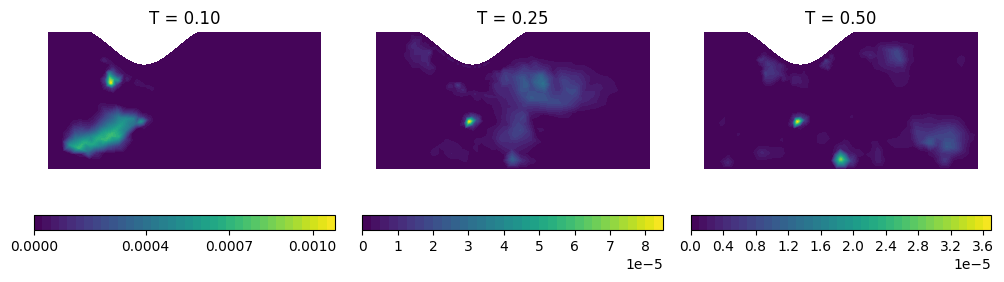}
  \caption{\textnormal{Test case 2, Advection-Diffusion problem in a 2D Stokes flow. Worst case scenario: bump close to the inflow ($\mu = 0.355$). First row: rollout prediction. Second row: RMSE related to each timestep between the prediction and the corresponding ground truth solution} }
  \label{fig:pred5bump}
\end{figure*}

%The best prediction, together with its ground truth solution, is shown in Figure \ref{fig:bestbump}.
\nicola{Indeed, as shown in Figures \ref{fig:pred12bump}-\ref{fig:pred5bump},} the predicted dynamics is still very accurate and we do not spot any propagation error. The prediction seems to get worse at some nodes which are either close to the bump or to the upper edge, in which we have imposed the no-slip conditions. \nicola{Conversely, the errors in proximity of the inflow are higher at initial times, but they tend to fade out as the simulation evolves (this is true also for our worst simulation, Figure \ref{fig:pred5bump}).}
%\begin{figure*}
% \centering
%  \includegraphics[scale=0.6]{Images/pred_2_bump2.png}
%  \includegraphics[scale=0.6]{Images/gt_2_bump2.png}
%  \caption{\textnormal{Test case 2, Advection-Diffusion problem in a 2D Stokes flow. Best model prediction ($\mu = 0.42$). First row: rollout prediction. Second row: ground truth solution.} }
%  \label{fig:bestbump}
%\end{figure*}
%The model performs well in predicting scenarios where the bump is located at new, unseen positions, as shown in Figure \ref{fig:pred12bump}. However, slightly higher errors can be spotted in nodes close to the inflow, despite they are kept under control and remain of the order $10^{-5}$. It is worth recalling that these initial errors tend to fade out, and the accuracy improves as the simulation evolves. This behavior is particularly evident when the bump is not located too close to the inflow, indeed higher errors are usually caused by particular positions of the bump.
%\pagebreak

%If however the position of the bump is too close to the inflow, it might heavily influence the prediction in the initial time steps. As a result, fixing the error in some nodes becomes increasingly difficult as time evolves. This can be seen quite well in Figure \ref{fig:pred5bump}, which represents a worst-case scenario for our predictions. The behavior of the nodes around the inflow has a significant impact on the accuracy of the simulation, as any error arising in this region can propagate throughout the domain. The proximity of the bump to the inflow also plays a crucial role in the self-adjustment of the simulation. However, the model overall performs very well, as demonstrated by the accuracy of its predictions.

\begin{figure}
\centering
  \includegraphics[width = 0.65\textwidth]{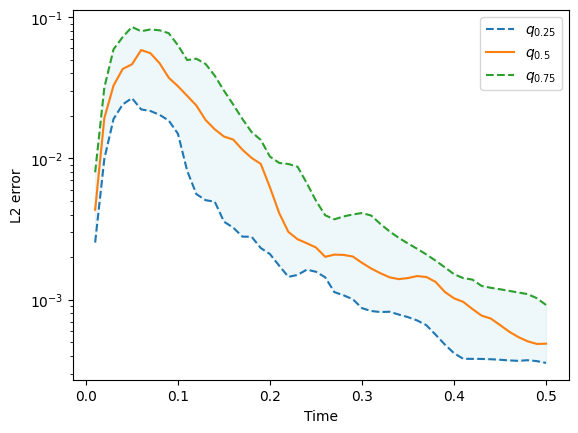}
  \caption{\textnormal{Test case 2, Advection-Diffusion problem in a 2D Stokes flow. $L^2$ relative error vs Time plot:  The dashed lines represent the first and the third quantiles of the $L^2$ errors among all the test predictions, while the orange line is the median. The shaded area can be seen as a confidence region for the simulation error.}}
  \label{fig:quantbump}
\end{figure}

Our qualitative considerations are also supported by the plot in Figure \ref{fig:quantbump}, which reports the behavior in time of the $L^2$ relative error between predictions and FOM solutions. \nicola{Clearly, the first time instants are the most challenging ones, as that is} when the inflow and the bump position determine the dynamics of the system. %; the third quantiles are instead influenced by the worst-case scenario.
Overall, we highlight that the model is able to self-adjust, since errors tend to decrease as the simulation time evolves, also showing some degree of robustness to the possible presence of noise during the simulation. 
%\par

%Furthermore, the decrease of the $L^2$ relative error in time indicates that the model is able to capture the system behavior more accurately at later time instants. This is expected since the initial condition and boundary conditions are better defined and the system behavior becomes more predictable as time evolves. Moreover, the third quantile in the plot shows that the model worst-case error is still relatively small, hence the model is robust to unexpected perturbations during the simulation. 

\nicola{Of note, these considerations hold uniformly over the test set, as clearly indicated by the width of the quantile bands.}
This is a desirable property since real-world problems often have some degree of uncertainty or noise, and a model that can handle %such
different scenarios %situations
is more likely to be useful in practice.

\subsubsection{Generalization to bumps with different positions and dimensions}
This example can be generalized by letting the bump \nicola{vary its dimension and possibly switch from the upper to the lower edge.} % also vary also along the lower edge, and its dimension change. 
Hence, we consider a new dataset consisting of 185 simulations in which the height of the bump is allowed to change, $h \in \{0.08,0.12,0.175\}$, and its center can vary along both the upper and lower edge in the interval $[0.4, 0.6]$, respectively. %In this way, we keep the bump sufficiently far from both the inflow and the outflow -- as seen before, this may entail some numerical errors in the GNN prediction. 
In other words, the new geometrical parameters are $\mub = (c_x,c_y,h) \in \Theta := [0.4,0.6] \times \{0.,0.5\} \times \{0.08,0.12,0.175\}$.  \\ % \par

%\begin{figure*}
% \centering
%  \includegraphics[scale=0.6]{Images/bump_ul_pred_30.png}
%  \includegraphics[scale=0.6]{Images/bump_ul_gt_30.png}
%  \caption{\textnormal{Test case 2, Advection-Diffusion problem in a 2D Stokes flow. Prediction obtained for $\mub = (0.528,0,0.12)$ with the bump in the lower edge with center at $x = 0.528$ and height $h = 0.12$. First row: rollout prediction. Second row: ground truth solution.}}
%  \label{fig:bumpul30}
%\end{figure*}

\begin{figure*}
 \centering
  \includegraphics[scale=0.6]{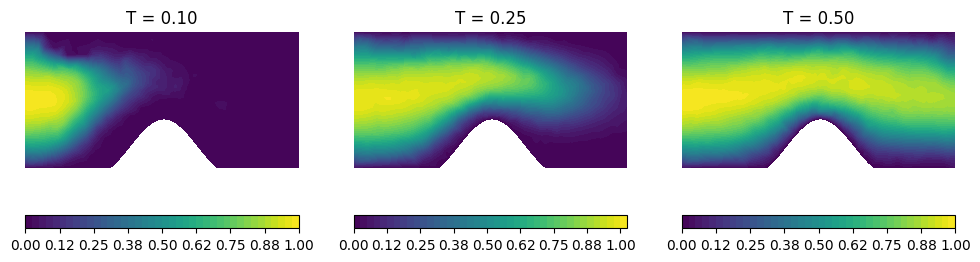}
  \includegraphics[scale=0.6]{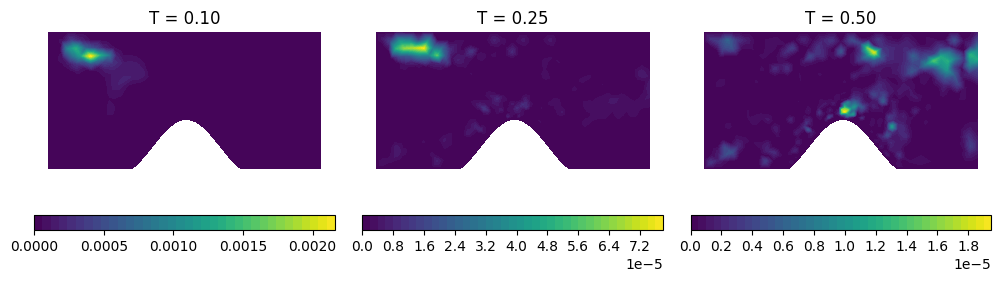}
  \caption{\textnormal{Test case 2, Advection-Diffusion problem in a 2D Stokes flow. Prediction obtained for $\mub = (0.453,0,0.175)$ with the bump in the lower edge with center at $x = 0.453$ and height $h = 0.175$. First row: rollout prediction. Second row: RMSE related to each timestep between the prediction and the corresponding ground truth solution.}}
  \label{fig:bumpul25}
\end{figure*}

\begin{figure*}
 \centering
  \includegraphics[scale=0.6]{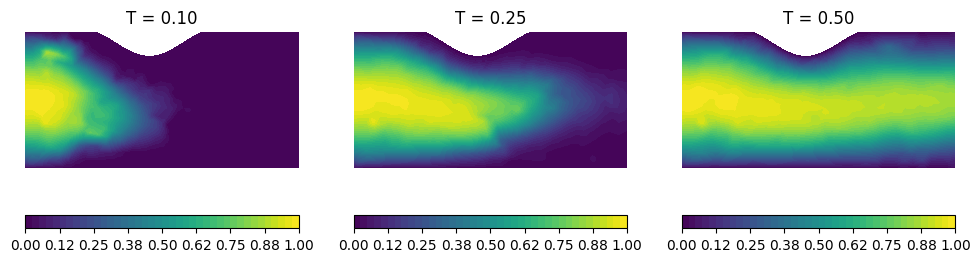}
  \includegraphics[scale=0.6]{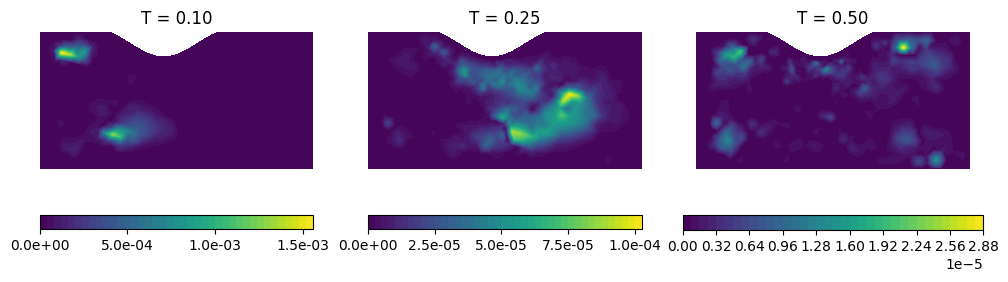}
  \caption{\textnormal{Test case 2, Advection-Diffusion problem in a 2D Stokes flow. Prediction with bump in the lower edge with center at $x = 0.453$ and height $h = 0.08$ ($\mub = (0.467,0.5,0.08)$). First row: rollout prediction. Second row: RMSE related to each timestep between the prediction and the corresponding ground truth solution.}}
  \label{fig:bumpul29}
\end{figure*}

\begin{figure}
\centering
  \includegraphics[width = 0.65\textwidth]{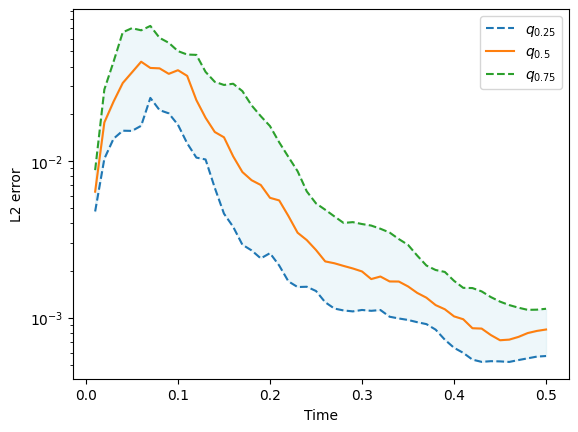}
  \caption{\textnormal{Test case 2, Advection-Diffusion problem in a 2D Stokes flow. $L^2$ relative error vs Time plot.}}
  \label{fig:quantbumpul}
\end{figure}

The results show that the implemented GNN-based model is able to learn correctly the geometry of the problem even if we the domain varies substantially within the dataset. %This is important since we can overcome overfitting, which is quite common when using deep neural networks. %In the following lines, some predictions of rollout simulations not seen during training are presented. In Figure \ref{fig:bumpul30}, the prediction computed by the GNN when the bump is on the lower edge with height $h = 0.12$ and center in $x = 0.528$, that is $\mub = (0.528,0,0.12)$, is reported. As in the previous example, the network shows some difficulty to predict the nodes close to the inflow but improves as the simulation time varies leading to a good approximation of the flow in the final steps. As we have previously seen, GNN-based models are often not good at predicting the regularity of the solution pattern, as we can see from the prediction at $T=0.1$ in Figure \ref{fig:bumpul30}. Here, we can also notice that the smoothness of the solution gets better as the simulation time increases. 
For instance, in Figure \ref{fig:bumpul25} the height of the bump influences a lot the system dynamics, however the network correctly infers the behavior of the flow around the obstacle. Here, the bump has height $h = 0.175$ and is located at the lower edge with center at $x = 0.453$, that is  $\mub = (0.453,0,0.175)$. The height of the bump has a significant impact on the accuracy of model prediction, particularly near the upper edge of the domain. Errors that arise in this region can propagate throughout the domain, affecting the accuracy of predictions at other locations as well. However, the self-adjustment mechanism of the model is effective in mitigating these errors as they propagate toward the outflow, resulting in improved accuracy in this region.
Overall, the model ability to account for the influence of the bump height on the flow dynamics contributes to its strong predictive performance.
In Figure \ref{fig:bumpul29} the bump has height $h = 0.08$ and is located at the upper edge with center at $x = 0.467$, that is $\mub = (0.467,0.5,0.08)$. The accuracy of the model predictions decreases when the size of the bump is smaller. This is primarily due to the fact that, as the size of the domain increases, so does the number of nodes, making the inference process more challenging. In this problem, the number of nodes varies from 936 to 1054, which is a wide range for unstructured meshes and geometries that differ significantly from each other. As a result, error propagation is more significant in this case compared to the other examples. This is highlighted by the persistence of relatively large errors at $T=0.25$, despite the overall dynamics being well-predicted by the model. This suggests that the model can effectively capture the underlying physics of the system, even in cases where the inference is more challenging due to the higher number of nodes. In general, the $L^{2}$-errors exhibit the same behavior as before: see Figure \ref{fig:quantbumpul} in comparison with Figure \ref{fig:quantbump}.
%Upon observing the $L^2$ relative error plot on the test set in Figure \ref{fig:quantbumpul}, we can draw some quantitative conclusions regarding the previously discussed results. The plot indicates that the test error has an appropriate upper bound and that it initially increases significantly during the first few time steps, which is consistent with the observed prediction behavior. We can see that the $L^2$ error trend is similar to that of Figure \ref{fig:quantbump}, even if it now decays more rapidly. This suggests that the self-adjustment property of the model is maintained even when the problem geometry exhibits greater variability.

\begin{figure*}
 \centering
  \includegraphics[scale=0.6]{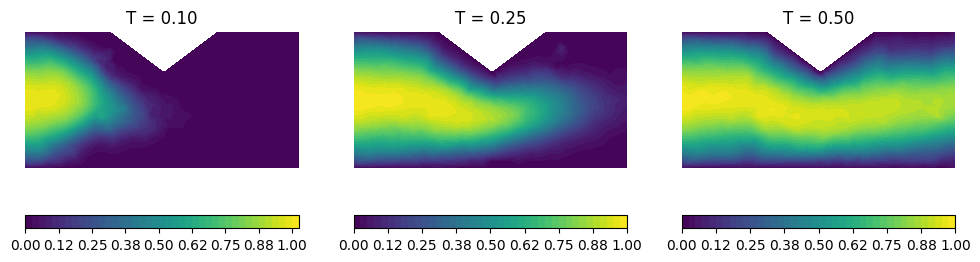}
  \includegraphics[scale=0.6]{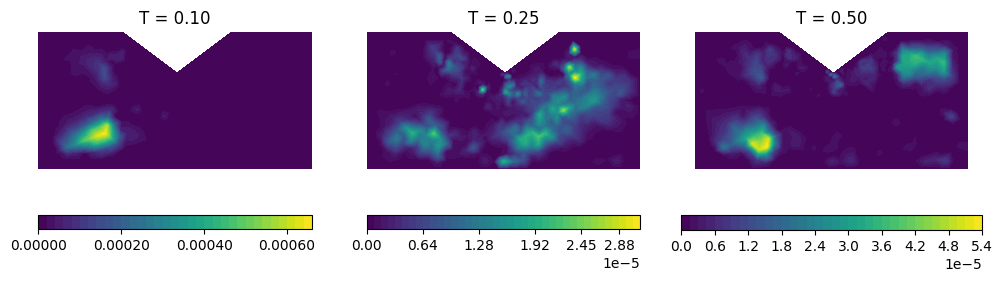}
  \caption{\textnormal{Test case 2, Advection-Diffusion problem in a 2D Stokes flow. Prediction with a triangular bump on the upper edge. First row: rollout prediction. Second row: RMSE related to each timestep between the prediction and the corresponding ground truth solution.}}
  \label{fig:bumptri}
\end{figure*}

We can further test the robustness of our model by evaluating its ability to predict solutions in channels with varying shapes of the bump, without requiring any retraining. Figure \ref{fig:bumptri} displays the prediction results of a simulation with a triangular bump located on the upper edge. In addition to the fact that the errors are of the same order of magnitude as previously discussed, the overall dynamics is, once again, accurately predicted. However, the regularity of the solution poses some difficulty for the model. Nonetheless, this does not appear to significantly impact the accuracy of the prediction.\par

This example highlights the flexibility of graph neural networks in handling simulations with variable geometries and limited training data,  while still producing reliable results.

\subsection{Advection-Diffusion problem in a 3D Stokes flow around a cylinder}
\label{sec:stokes3d}

To further increase the problem difficulty, we finally consider the same problem discussed in Section \ref{sec:stokes}, now set in a 3D domain obtained by an extrusion on the z-axis of the rectangle $R = (0,1) \times (0,0.5)$ with a cylindrical hole $$C = \{(x,y) : \ (x-c_{x})^2 + (y-c_{y})^2 \leq (0.05)^2\} .$$ We let the position of the obstacle vary as $\mub = (c_x,c_y) \in \Theta:=[0.2,0.4] \times [0.2,0.3]$. This time, we exploit the FOM to generate 150 different simulations\nicola{, 125 for training and 25 for testing.} %in the rectangle $R_c = [0.2,0.4] \times [0.2,0.3]$ resulting in 150 different simulations, that is . %An example of FOM solution can be seen in Figure \ref{fig:gtstokes3D}.
%\begin{figure*}
%\centering
%  \includegraphics[width=0.9\textwidth]{Images/stokes3d_gt.png}
%  \caption{\textnormal{Three different ground truth simulations.} }
%  \label{fig:gtstokes3D}
%\end{figure*} 

\subsubsection{Problem data}
The mesh nodes of the simulations vary from 1353 to 1542, thus increasing the complexity of the problem with respect to the other examples we have discussed so far. \nicola{We implement the approach following the same ideas adopted for the previous test cases: we refer to Table \ref{tab:hyperparameters} for further details about the network design and the training hyperparameters.} %except for the following ones: 

%and these are the other hyperparameters chosen: 

%\begin{itemize}
%\item $l = 32$, where $l$ is the latent dimension of the network;

%\item noise variance $\sigma^2 = 10^{-5}$;

%item total number of epochs $max\_epoch = 2000$;

%\item learning rate $\nu = 10^{-4}$ with decay $= 0.1$ after 500 and 1000 epochs;

%\item SiLU is chosen as the activation function for each MLP layer;

%\item number of MLP layers $mlp\_layers = 2$;

%\item message Passing steps $mp\_steps = 18$;

%\item loss weights $w_1 = 0.5, \ w_2 = 0.5$.
%\end{itemize} 

\nicola{In this regard, we have made some minor modifications to account for the increased complexity entailed by the presence of a three-dimensional geometry. These concern: an increased number of message-passing steps (to better cover the spatial domain), an increased noise variance (to further improve the stability of our simulations during rollout), and a reduced number of epochs (to avoid overfitting).} %In this third example, an additional computational burden hampers the complexity of the parameter-dependent dynamical system. To tackle these challenges, we have employed a specific training strategy where the model is trained for a relatively shorter period of 2000 epochs, with an increased number of message passing steps (18) to better capture the dynamics of the system. Additionally, we have also increased the noise variance from $10^{-6}$ to $10^{-5}$ to enable the model to learn handling higher levels of error propagation. This is particularly important as error propagation is a common issue when working with simulation rollouts for this problem. 
By adopting this training strategy, we aim to strike a balance between model accuracy and computational efficiency while still being able to capture the complex dynamics of the system.

\subsubsection{Numerical Results}
The results of the rollout predictions of the test simulations are summarized in Table \ref{table:errors}.
Clearly, the higher error obtained in this example is due to the increased complexity of the problem. However, it is noteworthy that despite the higher error, there is a significant improvement in time complexity. Once trained, our model can be up to two orders of magnitude faster  compared to the FOM solver. This reduction in time complexity can lead to faster and more efficient simulations, which is particularly important for time-critical applications or when a large number of simulations are required. Therefore, despite the slightly higher error, our model can still provide a significant advantage in terms of time and computational resources. %\par
%\begin{table}[H]
%    \caption*{\textbf{Advection-Diffusion problem in a 3D Stokes flow}}
%    \centering 
%    \begin{tabular}{|p{3em} c c c c c |}
%    \hline
    
\begin{figure*}
    % \centering
     \includegraphics[width=0.3\textwidth]{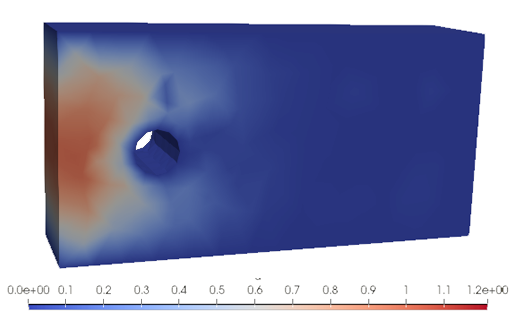}
     %\centering
     \includegraphics[width=0.3\textwidth]{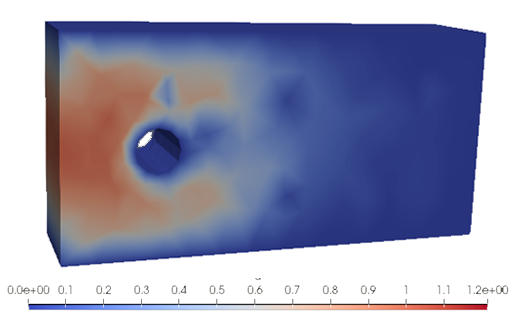}
     %\centering
     \includegraphics[width=0.3\textwidth]{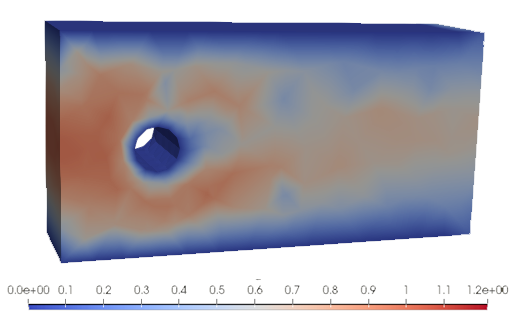}
    %\caption{Boxplots of the RMSEs of the two predictions}
    %\label{fig:stokes3d_36}
    % \centering
     \includegraphics[width=0.3\textwidth]{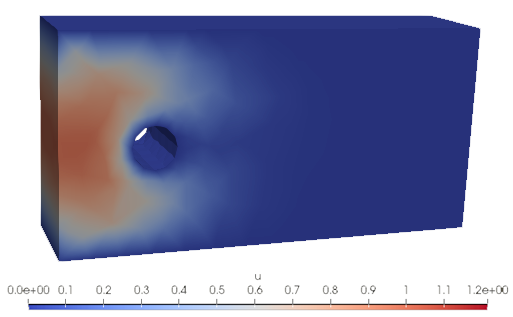}
     %\centering
     \;\;\;\;\;
     \includegraphics[width=0.3\textwidth]{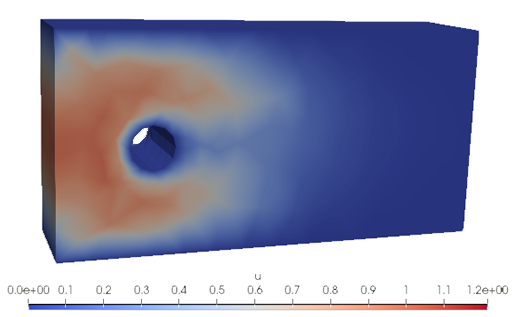}
     %\centering
     \;\;\;\;\;
     \includegraphics[width=0.3\textwidth]{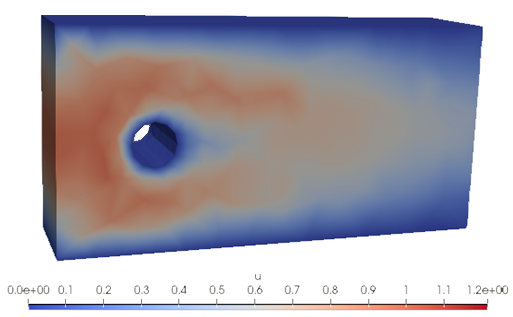}
    \caption{\textnormal{Test case 3, Advection-Diffusion problem in a 3D Stokes flow. Prediction obtained for $\mub = (0.29,0.25)$. 3 time steps of simulation. First row: Rollout prediction. Second row: ground truth solution}.}
    \label{fig:stokes3d_36}
\end{figure*}

%     & \textbf{RMSE (mean)} & \textbf{RMSE (max)} & \textbf{RMSE (min)} & $t_{\text{FOM}}$ & $t_{\text{GNN}}$  \\
%    \hline \hline
%    \textbf{AD 3} & $4.37 \mathrm{e}{-2}$ & $6.24 \mathrm{e}{-2}$ & $1.94 \mathrm{e}{-2}$ &$\approx 729.42$  & $\approx 10.4$  \\
%    
%    \hline
%    \end{tabular}
%    \\[10pt]
%    \caption{Test case 3, Advection-Diffusion problem in a 3D Stokes flow. Results of the test set predictions.}
%    \label{table:stokes3d}
%\end{table}
Upon examining the predictions in greater detail, as shown in Figure \ref{fig:stokes3d_36}, a comparison can be made between the prediction of a simulation with the obstacle positioned centrally, and its corresponding ground truth solution. It is evident that the simulation deteriorates as it progresses toward the outflow. Unlike the 2D case, self-adjustment is not observed in this scenario, as the nodes located to the right of the obstacle are heavily influenced by its position. This may result in some values being underestimated in the prediction, particularly in the tail of the flow. Unfortunately, this is a known drawback of \nicola{GNNs, %like the one implemented in this study, 
as deep architectures tend to oversmooth predictions.} % homogenize predictions if the network has too many nodes.
Therefore, even if the dynamics are predicted accurately, node values may be more dispersed. Furthermore, this problem is exacerbated by an increase in the number of message passing steps, which, in this example, are necessary for an acceptable prediction. \par

\begin{figure*}
     % \centering
     \includegraphics[width=0.3\textwidth]{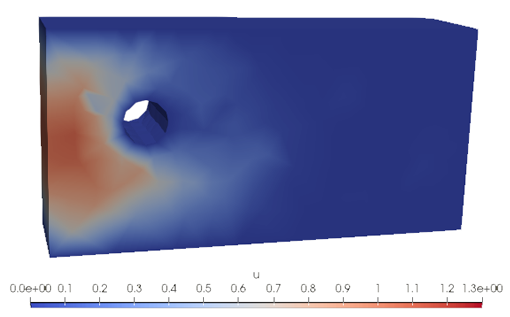}
     %\centering
     \includegraphics[width=0.3\textwidth]{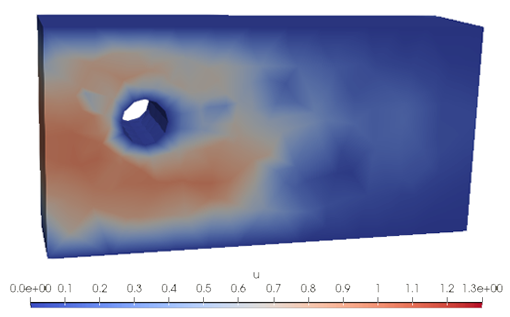}
     %\centering
     \includegraphics[width=0.3\textwidth]{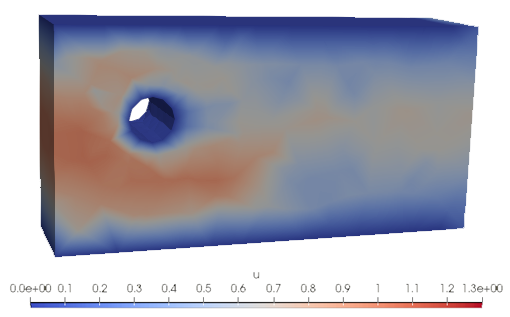}
    %\caption{Boxplots of the RMSEs of the two predictions}
    %\label{fig:stokes3d_36}
    % \centering
     \includegraphics[width=0.3\textwidth]{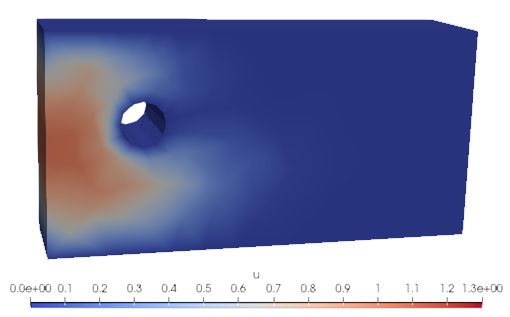}
     %\centering
     \;\;\;\;\;
     \includegraphics[width=0.3\textwidth]{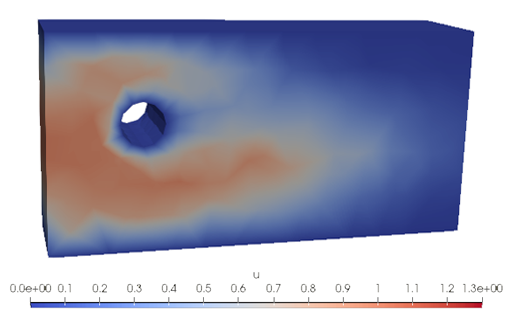}
     %\centering
     \;\;\;\;\;
     \includegraphics[width=0.3\textwidth]{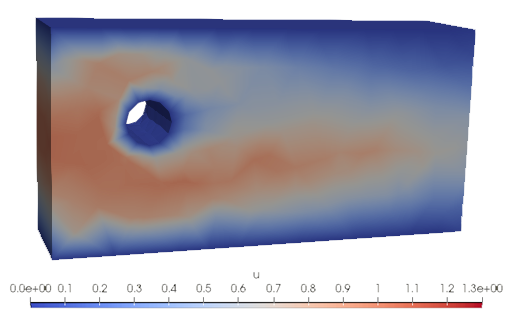}
    \caption{\textnormal{Test case 3, Advection-Diffusion problem in a 3D Stokes flow. Prediction obtained for $\mub = (0.23,0.3)$. 3 time steps of simulation. First row: Rollout prediction. Second row: ground truth solution}.}
    \label{fig:stokes3d_134}
\end{figure*}

Nevertheless, as illustrated in Figure \ref{fig:stokes3d_134}, modifying the position of the obstacle does not significantly affect the overall accuracy of the solution. Despite the aforementioned issues, the flow pattern is captured correctly, and no propagation of the errors is observed. Of remarkable importance is the consistently accurate prediction in the proximity of the obstacle, which is always a critical aspect to be predicted. This observation underscores the model ability to learn the geometrical properties of the problem while preserving the graph structure of the mesh. Therefore, these results suggest that the model is sufficiently robust in predicting flow patterns in various configurations, and can generalize well to other geometries. \par

\begin{figure}
\centering
  \includegraphics[width=0.65\textwidth]{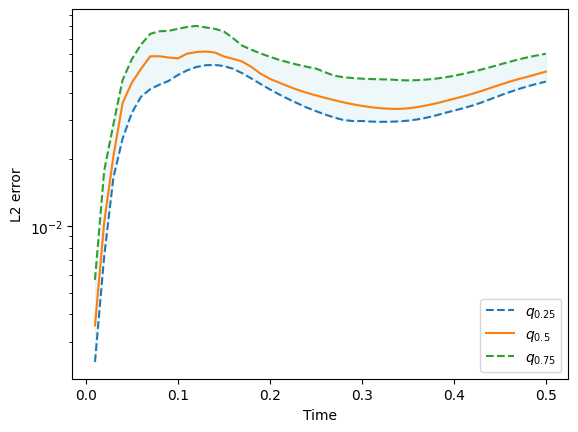}
  \caption{\textnormal{Test case 3, Advection-Diffusion problem in a 3D Stokes flow. $L^2$ relative error vs Time plot: The dashed lines represent the first and the third quantiles of the $L^2$ errors among all the test predictions, while the orange line is the median. The shaded area can be seen as a confidence region for the simulation error.}}
  \label{fig:quantstokes3d}
\end{figure}

Upon observing the $L^2$ relative error plot on the test set in Figure \ref{fig:quantstokes3d}, we can draw quantitative conclusions regarding the previously discussed results. The plot indicates that the test error has an appropriate upper bound and that it  increases significantly during the first few time steps, which is consistent with the observed prediction behavior. After the initial increase, the error gradually decays, showing that the model has learned the underlying dynamics of the system. However, towards the end of the simulation, we observe a slight increase in the error, which is coherent with what we have previously mentioned about the tendency of these architectures to dispersion. This behavior may be due to the accumulation of errors during the long-term prediction. Therefore, we can conclude that while the GNN-based model shows promising results, there is still room for improvement in terms of accuracy and robustness.

\subsection{Comparison with Feed Forward Neural Networks}
Feed Forward Neural Networks (FFNNs) are usually employed for building reduced-order models because they have the capability to capture strong nonlinearity through their fully connected structure \cite{hesthaven2018non, mucke2021reduced, kelly2022data}. \nicola{As we explained in the introduction, however, using FFNNs is not straightforward when dealing with geometric variability, as these models require fixing both the input and output dimension. We recall, in fact, that FFNN architectures are nothing but MLP units. More precisely, any FFNN model of depth $s\ge1$ is a map of the form
\begin{equation}
    \label{eq:mlp}
    \Psi:= L_{s+1}\circ L_{s}\circ \dots \circ L_{1},
\end{equation}
where $L_{i}:\mathbb{R}^{n_{i}}\to\mathbb{R}^{n_{i+1}}$ are nonlinear maps (layers) operating as
$$L_{i}:\mathbf{v}\mapsto\rho_{i}\left(\mathbf{W}_{i}\mathbf{v}+\mathbf{b}_{i}\right),$$
where $\rho_{i}:\mathbb{R}\to\mathbb{R}$ is the activation function (acting componentwise), while $\mathbf{W}_{i}$ and $\mathbf{b}_{i}$ are the trainable parameters (weights and biases, respectively). Compared to GNNs, these computational units are substantially less flexible, as they can only accept inputs of a specific dimension (here, $n_{1}$), and will always produce outputs of a given size (here, $n_{s+2}$). Consequently, it would be impossible to implement a model such as \eqref{eq:time-stepping} using naive FFNNs, at least not in the case of parametrized domains. 
}

\nicola{Still, it is true that we could circumvent this problem by interpolating each PDE solution over a fixed rectangular grid, with the convention that $u_{\mub}(\mathbf{x})=0$ if $\mathbf{x}\notin\Omega_{\mub}$. Then, it would be possible -in principle- to replicate the same construction proposed in Section \ref{sec:surrogate}, but only using dense architectures such as \eqref{eq:mlp}. Our claim, however, is that without the \textit{relational inductive bias} of GNNs, these model have no hope of generalizing over unseen geometries. To prove this, we shall compare the performances obtained by using either FFNNs or GNNs in \eqref{eq:time-stepping}}
%In this final section we compare the performances \nicola{of classical FFNNs } of common FFNN to our GNN-based model, 
on the two examples of Sections \ref{sec:ad} and \ref{sec:stokes}. In particular, to simplify, we consider for each problem the following datasets: 
\begin{itemize}
    \item for the first Advection-Diffusion problem, we let the obstacle vary only in its position, resulting in 100 random simulations (80 for training and 20 for testing) ;
    
    \item for the Advection-Diffusion problem in a 2D Stokes flow, we let the bump vary in its position along the upper and lower edges but not in its height, resulting in 125 random simulations (100 for training and 25 for testing).
\end{itemize} \par

\begin{figure*}
      \centering
\includegraphics[width=0.495\textwidth]{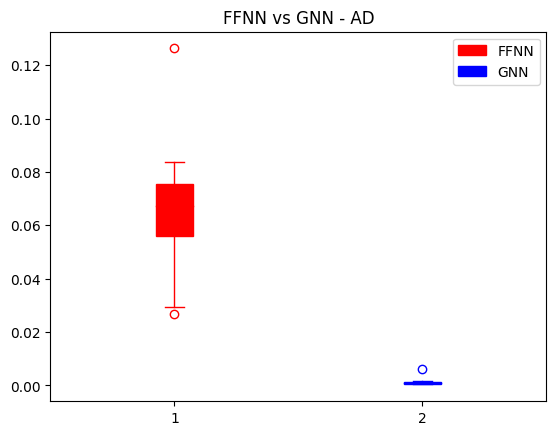}    
\includegraphics[width=0.495\textwidth]{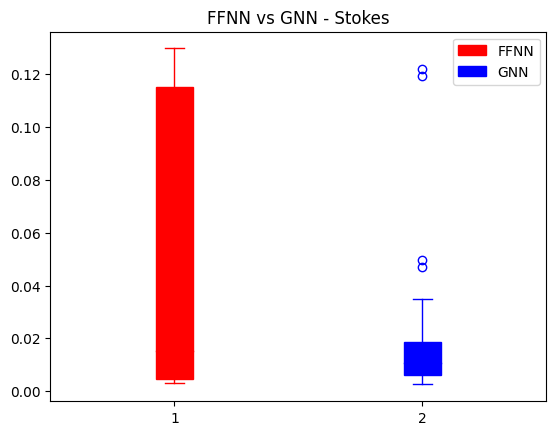} 
       \caption{Boxplots of the RMSEs of the two predictions}
    \label{fig:boxplotffnn}
\end{figure*}

\begin{figure*}
 \centering
  \includegraphics[scale=0.6]{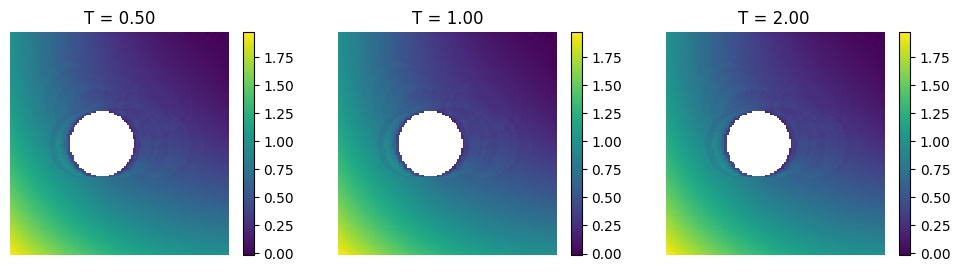}
  \includegraphics[scale=0.6]{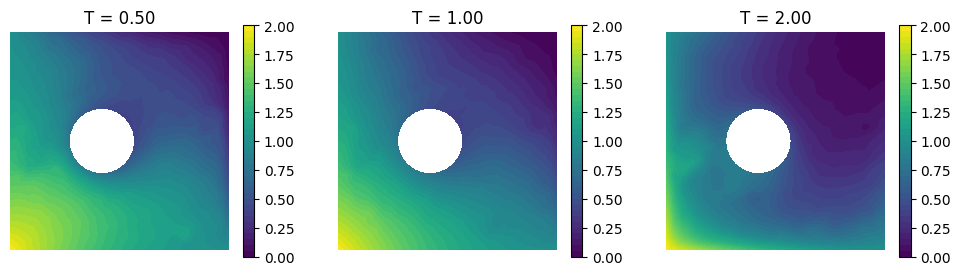}
  \caption{\textnormal{Test case 1, Advection - Diffusion problem. Comparison between FFNN and GNN prediction for $\mub = (c_x,c_y) = (0.4,0.5)$. First row: FFNN prediction. Second row: GNN prediction.} }
  \label{fig:adcomp}
\end{figure*}

\begin{figure*}
 \centering
  \includegraphics[scale=0.37]{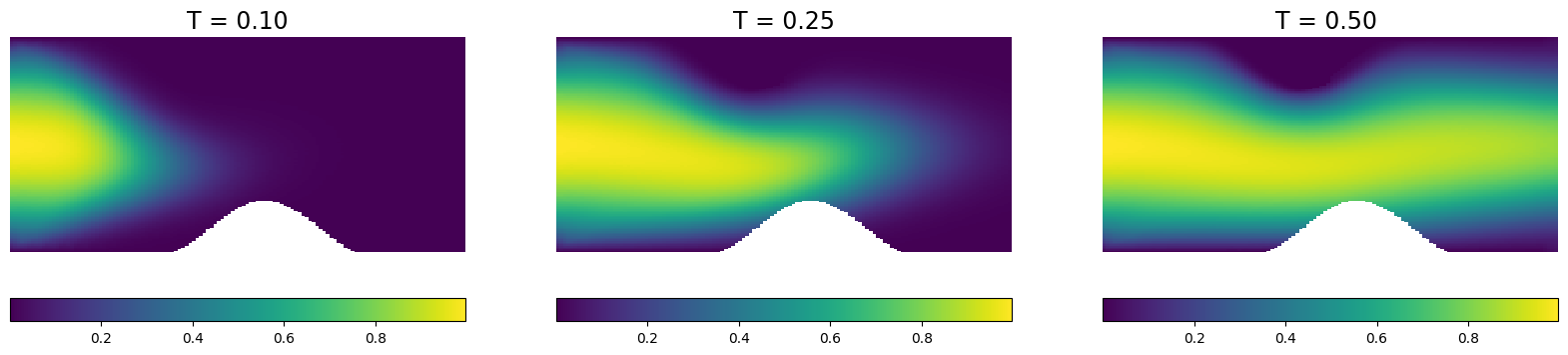}
  \includegraphics[scale=0.6]{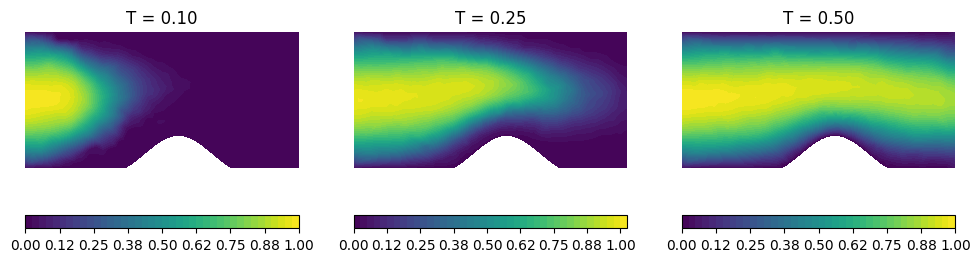}
  \caption{\textnormal{Test case 2, Advection - Diffusion problem in a 2D Stokes flow. Comparison between FFNN and GNN prediction for $\mub = (c_x,c_y) = (0.6,0)$. First row: FFNN prediction. Second row: GNN prediction}.}
  \label{fig:bumpcomp}
\end{figure*}

 %In order to train the FFNN, we need to make the size of all the data compatible, i.e., we need to map all the simulation data onto a common set of degrees of freedom. Therefore, we interpolate the region of interest onto a modeling $128 \times 128$ vertices grid so that the order of the nodes in the mesh is preserved.
 In order to train the FFNNs, we map all the simulations onto a common rectangular grid consisting of $128 \times 128$ vertices. We then repeat the same construction presented in Section \ref{sec:surrogate}, up to replacing the GNN modules with MLP layers.
%Then, we build an Encoder-Decoder fully connected network to approximate the solution of the system at time $t^{n+1}$, given the one at time $t^n$, as described in Equation (\ref{eq:gnn_vec}).
%The training is performed by one-step prediction and minimizing the batch loss in  (\ref{eq:loss}), as already presented in Section \ref{sec:train_alg}. Moreover, we pass as input to the network a $N_{batch} \times 128 \times 128 \times 3$ tensor containing the same features as the GNN model. 
Since FFNNs tend to overfit if trained for a long time, we train the model for $500$ epochs, using the same learning rate and loss weights as the ones described in Sections \ref{sec:ad} and \ref{sec:stokes}.
\\\\
\nicola{Results are in Figures \ref{fig:boxplotffnn}-\ref{fig:bumpcomp}.}
%The prediction at testing is done by exploiting the rollout of the simulation as shown in Equation (\ref{eq:pred}) in order to compare the robustness to propagation errors of the two models. 
\nicola{As testified by the boxplots, the quality of FFNN predictions can change a lot from case to case; conversely, GNN-surrogates are much more stable and report errors with a smaller variability.} %In Figure \ref{fig:boxplotffnn} we can clearly see that the predictions done with GNN show less variance than the ones obtained  with the FFNN. Even if the range of the errors is similar, the GNN errors are overall better. % despite some anomalies. 
Two examples of FFNN predictions on test set simulations, together with the corresponding GNN predictions, are shown in Figures \ref{fig:adcomp} and \ref{fig:bumpcomp}. The predictions are done using the same values of geometrical parameters in order to highlight the different performances in the generalization on unseen domains. %A significant difference between FFNN and GNN models for solving problems that require making inference on the geometry clearly arises.
\nicola{Here, the difference between the two approaches becomes evident.} While FFNNs may capture the overall dynamics of the system quite well, they \nicola{fail at understanding the geometrical properties of the solution, which ultimately makes them unable to generalize (see, e.g., Figure \ref{fig:bumpcomp}, where the FFNN clearly ignores the actual location of the bump).} %are not able to handle any geometric property of the solution. 
%This limitation is particularly evident in simple examples, therefore implying that the benefits of using a GNN can be even more pronounced when dealing with more complex problems.

Another key aspect to analyze is the number of model parameters. In fact, due to their fully connected structure, the complexity of FFNNs can increase dramatically with the problem complexity, implying a higher tendency to overfitting (thus the need of more training data) and less scalable models. 

A possible strategy to overcome both these issues could be to rely on grid-based models, such as Convolutional Neural Networks (CNNs). In fact, these models can reduce the number of parameters by sharing them, which helps to mitigate the overfitting issue. %In Equation \ref{eq:mlp}, the weight matrix $\mathbf{W}$ and bias vector $\mathbf{b}$ are shared across simulations in the CNN architecture. However, this may lead to lower prediction accuracy due to excessive generalization. This is also a drawback of GNN models, which often struggle to propagate information globally across the mesh. %\par
However, while CNNs might be able to capture the dynamics of the system, they would still ignore the geometric structure of the problem, which makes them unsuitable for complex geometries. %Therefore, while CNNs are a viable alternative to FFNNs, they are not a perfect solution for problems that require the consideration of both geometry and dynamics, such as the ones we have discussed.
Similarly, alternative approaches such as %In this context, 
Mesh-Informed Neural Networks (MINNs) \cite{franco2022learning} do not provide a comprehensive solution, as they can only tackle one geometry at a time. %In fact, although MINNs can incorporate geometric information into their definition, they require the underlying spatial domain to have a fixed shape. %be fixed cannot handle However, they are still limited in their ability to interpolate the solution in a more general domain, where the number of nodes in the mesh may vary.

%% file: aipsamp.bbl
\begin{thebibliography}{10}

\bibitem{amsallem2011online}
David Amsallem and Charbel Farhat.
\newblock An online method for interpolating linear parametric reduced-order
  models.
\newblock {\em SIAM Journal on Scientific Computing}, 33(5):2169--2198, 2011.

\bibitem{10.5555/3157382.3157601}
Peter Battaglia, Razvan Pascanu, Matthew Lai, Danilo~Jimenez Rezende, and Koray
  kavukcuoglu.
\newblock Interaction networks for learning about objects, relations and
  physics.
\newblock In {\em Proceedings of the 30th International Conference on Neural
  Information Processing Systems}, NIPS'16, page 4509–4517, Red Hook, NY,
  USA, 2016. Curran Associates Inc.

\bibitem{battaglia2018relational}
Peter~W Battaglia, Jessica~B Hamrick, Victor Bapst, Alvaro Sanchez-Gonzalez,
  Vinicius Zambaldi, Mateusz Malinowski, Andrea Tacchetti, David Raposo, Adam
  Santoro, Ryan Faulkner, et~al.
\newblock Relational inductive biases, deep learning, and graph networks.
\newblock {\em arXiv preprint arXiv:1806.01261}, 2018.

\bibitem{bonito2021nonlinear}
Andrea Bonito, Albert Cohen, Ronald DeVore, Diane Guignard, Peter Jantsch, and
  Guergana Petrova.
\newblock Nonlinear methods for model reduction.
\newblock {\em ESAIM: Mathematical Modelling and Numerical Analysis},
  55(2):507--531, 2021.

\bibitem{brivio2023error}
Simone Brivio, Stefania Fresca, Nicola~Rares Franco, and Andrea Manzoni.
\newblock Error estimates for pod-dl-roms: a deep learning framework for
  reduced order modeling of nonlinear parametrized pdes enhanced by proper
  orthogonal decomposition.
\newblock {\em arXiv preprint arXiv:2305.04680}, 2023.

\bibitem{buza2021using}
Gergely Buza, Shobhit Jain, and George Haller.
\newblock Using spectral submanifolds for optimal mode selection in nonlinear
  model reduction.
\newblock {\em Proceedings of the Royal Society A}, 477(2246):20200725, 2021.

\bibitem{carlberg2015adaptive}
Kevin Carlberg.
\newblock Adaptive h-refinement for reduced-order models.
\newblock {\em International Journal for Numerical Methods in Engineering},
  102(5):1192--1210, 2015.

\bibitem{cicci2022deep}
Ludovica Cicci, Stefania Fresca, and Andrea Manzoni.
\newblock Deep-hyromnet: A deep learning-based operator approximation for
  hyper-reduction of nonlinear parametrized pdes.
\newblock {\em Journal of Scientific Computing}, 93(2):57, 2022.

\bibitem{fatone2022long}
Federico Fatone, Stefania Fresca, and Andrea Manzoni.
\newblock Long-time prediction of nonlinear parametrized dynamical systems by
  deep learning-based reduced order models.
\newblock {\em arXiv preprint arXiv:2201.10215}, 2022.

\bibitem{franco2023deep}
Nicola Franco, Andrea Manzoni, and Paolo Zunino.
\newblock A deep learning approach to reduced order modelling of parameter
  dependent partial differential equations.
\newblock {\em Mathematics of Computation}, 92(340):483--524, 2023.

\bibitem{franco2022learning}
Nicola~Rares Franco, Andrea Manzoni, and Paolo Zunino.
\newblock Learning operators with mesh-informed neural networks.
\newblock {\em arXiv preprint arXiv:2203.11648}, 2022.

\bibitem{fresca2021comprehensive}
Stefania Fresca, Luca Dede’, and Andrea Manzoni.
\newblock A comprehensive deep learning-based approach to reduced order
  modeling of nonlinear time-dependent parametrized pdes.
\newblock {\em Journal of Scientific Computing}, 87:1--36, 2021.

\bibitem{fresca2022pod}
Stefania Fresca and Andrea Manzoni.
\newblock Pod-dl-rom: enhancing deep learning-based reduced order models for
  nonlinear parametrized pdes by proper orthogonal decomposition.
\newblock {\em Computer Methods in Applied Mechanics and Engineering},
  388:114181, 2022.

\bibitem{gladstone2023gnn}
Rini~Jasmine Gladstone, Helia Rahmani, Vishvas Suryakumar, Hadi Meidani, Marta
  D'Elia, and Ahmad Zareei.
\newblock Gnn-based physics solver for time-independent pdes.
\newblock {\em arXiv preprint arXiv:2303.15681}, 2023.

\bibitem{guo2019data}
Mengwu Guo and Jan~S Hesthaven.
\newblock Data-driven reduced order modeling for time-dependent problems.
\newblock {\em Computer methods in applied mechanics and engineering},
  345:75--99, 2019.

\bibitem{guo2022bayesian}
Mengwu Guo, Shane~A McQuarrie, and Karen~E Willcox.
\newblock Bayesian operator inference for data-driven reduced-order modeling.
\newblock {\em Computer Methods in Applied Mechanics and Engineering},
  402:115336, 2022.

\bibitem{hamilton2017inductive}
W.~L. Hamilton, R.~Ying, and J.~Leskovec.
\newblock Inductive representation learning on large graphs.
\newblock {\em Advances in Neural Information Processing Systems}, pages
  1024--1024, 2017.

\bibitem{HamiltonYL17}
W.~L. Hamilton, R.~Ying, and J.~Leskovec.
\newblock Inductive representation learning on large graphs.
\newblock {\em CoRR}, abs/1706.02216, 2017.

\bibitem{hess2019localized}
Martin Hess, Alessandro Alla, Annalisa Quaini, Gianluigi Rozza, and Max
  Gunzburger.
\newblock A localized reduced-order modeling approach for pdes with bifurcating
  solutions.
\newblock {\em Computer Methods in Applied Mechanics and Engineering},
  351:379--403, 2019.

\bibitem{hesthaven2022rank}
Jan~S Hesthaven, Cecilia Pagliantini, and Nicol{\`o} Ripamonti.
\newblock Rank-adaptive structure-preserving model order reduction of
  hamiltonian systems.
\newblock {\em ESAIM: Mathematical Modelling and Numerical Analysis},
  56(2):617--650, 2022.

\bibitem{hesthaven2016certified}
Jan~S Hesthaven, Gianluigi Rozza, Benjamin Stamm, et~al.
\newblock {\em Certified reduced basis methods for parametrized partial
  differential equations}, volume 590.
\newblock Springer, 2016.

\bibitem{hesthaven2018non}
Jan~S Hesthaven and Stefano Ubbiali.
\newblock Non-intrusive reduced order modeling of nonlinear problems using
  neural networks.
\newblock {\em Journal of Computational Physics}, 363:55--78, 2018.

\bibitem{horie2022physics}
Masanobu Horie and Naoto Mitsume.
\newblock Physics-embedded neural networks: Graph neural pde solvers with mixed
  boundary conditions.
\newblock {\em Advances in Neural Information Processing Systems},
  35:23218--23229, 2022.

\bibitem{kelly2022data}
Sean~T Kelly and Bogdan~I Epureanu.
\newblock Data-driven reduced-order model for turbomachinery blisks with
  friction nonlinearity.
\newblock In {\em Nonlinear Structures \& Systems, Volume 1: Proceedings of the
  40th IMAC, A Conference and Exposition on Structural Dynamics 2022}, pages
  97--100. Springer, 2022.

\bibitem{kingma2015adam}
D.~P. Kingma and J.~Ba.
\newblock {\em Adam: A method for stochastic optimization}.
\newblock International Conference on Learning Representations (ICLR), 2015.

\bibitem{kipf2016semi}
T.N. Kipf and M.~Welling.
\newblock Semi-supervised classification with graph convolutional networks.
\newblock {\em 5th International Conference on Learning Representations
  (ICLR-17)}, 2016.

\bibitem{li2023model}
Mingwu Li, Shobhit Jain, and George Haller.
\newblock Model reduction for constrained mechanical systems via spectral
  submanifolds.
\newblock {\em Nonlinear Dynamics}, 111(10):8881--8911, 2023.

\bibitem{mucke2021reduced}
Nikolaj~T M{\"u}cke, Sander~M Boht{\'e}, and Cornelis~W Oosterlee.
\newblock Reduced order modeling for parameterized time-dependent pdes using
  spatially and memory aware deep learning.
\newblock {\em Journal of Computational Science}, 53:101408, 2021.

\bibitem{murtagh1991multilayer}
Fionn Murtagh.
\newblock Multilayer perceptrons for classification and regression.
\newblock {\em Neurocomputing}, 2(5-6):183--197, 1991.

\bibitem{negri2013reduced}
Federico Negri, Gianluigi Rozza, Andrea Manzoni, and Alfio Quarteroni.
\newblock Reduced basis method for parametrized elliptic optimal control
  problems.
\newblock {\em SIAM Journal on Scientific Computing}, 35(5):A2316--A2340, 2013.

\bibitem{pmlr-v48-niepert16}
Mathias Niepert, Mohamed Ahmed, and Konstantin Kutzkov.
\newblock Learning convolutional neural networks for graphs.
\newblock In Maria~Florina Balcan and Kilian~Q. Weinberger, editors, {\em
  Proceedings of The 33rd International Conference on Machine Learning},
  volume~48 of {\em Proceedings of Machine Learning Research}, pages
  2014--2023, New York, New York, USA, 20--22 Jun 2016. PMLR.

\bibitem{pagliantini2021dynamical}
Cecilia Pagliantini.
\newblock Dynamical reduced basis methods for hamiltonian systems.
\newblock {\em Numerische Mathematik}, 148(2):409--448, 2021.

\bibitem{pegolotti2023learning}
Luca Pegolotti, Martin~R Pfaller, Natalia~L Rubio, Ke~Ding, Rita~Brugarolas
  Brufau, Eric Darve, and Alison~L Marsden.
\newblock Learning reduced-order models for cardiovascular simulations with
  graph neural networks.
\newblock {\em arXiv preprint arXiv:2303.07310}, 2023.

\bibitem{pfaff2020learning}
Tobias Pfaff, Meire Fortunato, Alvaro Sanchez-Gonzalez, and Peter~W Battaglia.
\newblock Learning mesh-based simulation with graph networks.
\newblock {\em arXiv preprint arXiv:2010.03409}, 2020.

\bibitem{pichi2023graph}
Federico Pichi, Beatriz Moya, and Jan~S Hesthaven.
\newblock A graph convolutional autoencoder approach to model order reduction
  for parametrized pdes.
\newblock {\em arXiv preprint arXiv:2305.08573}, 2023.

\bibitem{quarteroni1994numerical}
A.~Quarteroni and A.~Valli.
\newblock {\em Numerical approximation of partial differential equations},
  volume~23.
\newblock Springer, 01 1994.

\bibitem{quarteroni2015reduced}
Alfio Quarteroni, Andrea Manzoni, and Federico Negri.
\newblock {\em Reduced basis methods for partial differential equations: an
  introduction}, volume~92.
\newblock Springer, 2015.

\bibitem{romor2023non}
Francesco Romor, Giovanni Stabile, and Gianluigi Rozza.
\newblock Non-linear manifold reduced-order models with convolutional
  autoencoders and reduced over-collocation method.
\newblock {\em Journal of Scientific Computing}, 94(3):74, 2023.

\bibitem{rumelhart2986learning}
D.~Rumelhart, G.~Hinton, and R.~Williams.
\newblock Learning representations by back-propagating errors.
\newblock {\em Nature}, pages 533--536, 1986.

\bibitem{rusch2023survey}
T~Konstantin Rusch, Michael~M Bronstein, and Siddhartha Mishra.
\newblock A survey on oversmoothing in graph neural networks.
\newblock {\em arXiv preprint arXiv:2303.10993}, 2023.

\bibitem{sanchez2020learning}
Alvaro Sanchez-Gonzalez, Jonathan Godwin, Tobias Pfaff, Rex Ying, Jure
  Leskovec, and Peter Battaglia.
\newblock Learning to simulate complex physics with graph networks.
\newblock In {\em International conference on machine learning}, pages
  8459--8468. PMLR, 2020.

\bibitem{scarselli2008graph}
Franco Scarselli, Marco Gori, Ah~Chung Tsoi, Markus Hagenbuchner, and Gabriele
  Monfardini.
\newblock The graph neural network model.
\newblock {\em IEEE transactions on neural networks}, 20(1):61--80, 2008.

\bibitem{ijcai2021p214}
Yunsheng Shi, Zhengjie Huang, Shikun Feng, Hui Zhong, Wenjing Wang, and Yu~Sun.
\newblock Masked label prediction: Unified message passing model for
  semi-supervised classification.
\newblock In Zhi-Hua Zhou, editor, {\em Proceedings of the Thirtieth
  International Joint Conference on Artificial Intelligence, {IJCAI-21}}, pages
  1548--1554. International Joint Conferences on Artificial Intelligence
  Organization, 8 2021.
\newblock Main Track.

\bibitem{shlomi2020graph}
Jonathan Shlomi, Peter Battaglia, and Jean-Roch Vlimant.
\newblock Graph neural networks in particle physics.
\newblock {\em Machine Learning: Science and Technology}, 2(2):021001, 2020.

\bibitem{shukla2022scalable}
Khemraj Shukla, Mengjia Xu, Nathaniel Trask, and George~E Karniadakis.
\newblock Scalable algorithms for physics-informed neural and graph networks.
\newblock {\em Data-Centric Engineering}, 3:e24, 2022.

\bibitem{veli2017graph}
P.~Veličković, G.~Cucurull, A.~Casanova, A.~Romero, P.~Liò, and Y.~Bengio.
\newblock Graph attention networks.
\newblock {\em 6th International Conference on Learning Representation}, 2017.

\bibitem{zhang2019graph}
Si~Zhang, Hanghang Tong, Jiejun Xu, and Ross Maciejewski.
\newblock Graph convolutional networks: a comprehensive review.
\newblock {\em Computational Social Networks}, 6(1):1--23, 2019.

\end{thebibliography}
